\documentclass[11pt]{article}
\usepackage{latexsym}
\usepackage{amsmath}
\usepackage{amssymb}
\usepackage{amsfonts}
\usepackage{eucal}
\newtheorem{Theorem}{Theorem}[section]
\newtheorem{Proposition}[Theorem]{Proposition}

\newtheorem{Corollary}[Theorem]{Corollary}
\newtheorem{Lemma}[Theorem]{Lemma}
\newtheorem{Remark}[Theorem]{Remark}

\def\R{\mathbb{R}}
\def\D{\mathcal{D}}
\def\EE{\mathrm{E}}
\def\X{\bold{X}}
\begin{document}
\title{Large deviations of Poisson cluster processes
}
\author{Charles Bordenave\thanks{INRIA/ENS, D{\'e}partment d'Informatique,
45 rue d'Ulm, F-75230 Paris Cedex 05, France. e-mail:
\tt{charles.bordenave@ens.fr}}\,\, and \,\,Giovanni Luca
Torrisi\thanks{Istituto per le Applicazioni del Calcolo "Mauro
Picone" (IAC), Consiglio Nazionale delle Ricerche (CNR), Viale del
Policlinico 137, I-00161 Roma, Italia. e-mail:
\tt{torrisi@iac.rm.cnr.it}}}
\date{}
\maketitle

\begin{abstract}
In this paper we prove scalar and sample path large deviation
principles for a large class of Poisson cluster processes.~As a
consequence, we provide a large deviation principle for ergodic
Hawkes point processes.
\end{abstract}


\section{Introduction}
\label{sec:int} Poisson cluster processes are one of the most
important classes of point process models (see Daley and
Vere-Jones (2003) and M{\o}ller and Waagepetersen (2004)).~They
are natural models for the location of objects in the space, and
are widely used in point process studies whether theoretical or
applied.~Very popular and versatile Poisson cluster processes are
the so-called self-exciting or Hawkes processes (Hawkes (1971a),
(1971b); Hawkes and Oakes (1974)).~From a theoretical point of
view Hawkes processes combine both a Poisson cluster process
representation and a simple stochastic intensity representation.

Poisson cluster processes found applications in cosmology, ecology
and epidemiology; see, respectively, Neyman and Scott (1958), Brix
and Chadoeuf (2002) and M{\o}ller (2003).~Hawkes processes are
particularly appealing for seismological applications.~Indeed,
they are widely used as statistical models for the standard
activity of earthquake series; see the papers by Ogata and Akaike
(1982), Vere-Jones and Ozaki (1982), Ogata (1988) and Ogata
(1998).~Hawkes processes have also aspects appealing to
neuroscience applications; see the paper by Johnson (1996).~More
recently, Hawkes processes found applications to finance, see
Chavez-Demoulin, Davison and Mc Neil (2005), and to DNA modeling,
see Gusto and Schbath (2005).

In this paper we derive scalar and sample path large deviation
principles for Poisson cluster processes.~The paper is organized
as follows.~In Section 2 we give some preliminaries on Poisson
cluster processes, Hawkes processes and large deviations.~In
Section 3 we provide scalar large deviation principles for Poisson
cluster processes, under a light-tailed assumption on the number
of points per cluster.~As consequence, we provide scalar large
deviations for ergodic Hawkes processes.~Section 4 is devoted to
sample path large deviations of Poisson cluster processes.~First,
we prove a sample path large deviation principle on $D[0,1]$
equipped with the topology of point-wise convergence, under a
light-tailed assumption on the number of points per
cluster.~Second, we give a sample path large deviation principle
on $D[0,1]$ equipped with the topology of uniform convergence,
under a super-exponential assumption on the number of points per
cluster.~In Section 5 we prove large deviations for spatial
Poisson cluster processes, and we provide the asymptotic behavior
of the void probability function and the empty space function.~We
conclude the paper with a short discussion.

\section{Preliminaries}
\label{sec:prel} In this section we recall the definition of
Poisson cluster process, Hawkes process, and the notion of large
deviation principle.

\subsection{Poisson cluster processes}
A Poisson cluster process $\bold{X}\subset\mathbb{R}$ is a point
process.~The clusters centers of $\bold{X}$ are given by
particular points called immigrants; the other points of the
process are called offspring.~The formal definition of the process
is the following:

\begin{itemize}
\item[(a)] The immigrants are distributed according to a
homogeneous Poisson process $I$ with points $X_i\in\mathbb{R}$ and
intensity $\nu>0$.

\item[(b)] Each immigrant $X_i$ generates a cluster $C_i=C_{X_i}$,
which is a finite point process containing $X_i$.

\item[(c)] Given the immigrants, the centered clusters

\[
C_i-X_i=\{Y-X_i:\,Y\in C_i\},\quad\text{$X_i\in I$}
\]

\noindent are independent, identically distributed (iid for
short), and independent of $I$.

\item[(d)] $\bold{X}$ consists of the union of all clusters.
\end{itemize}

\noindent The number of points in a cluster is denoted by $S$.~We
will assume that $\EE[S]<\infty$.~Let $\bold{Y}$ be a point
process on $\mathbb{R}$ and $N_{\bold{Y}}(0,t]$ the number of
points of $\bold{Y}$ in the interval $(0,t]$.~$\bold{Y}$ is said
stationary if its law is translations invariant, is said ergodic
if it is stationary, with a finite intensity
$\mathrm{E}[N_{\bold{Y}}(0,1]]$, and

\[
\lim_{t\rightarrow\infty}\frac{N_{\bold{Y}}(0,t]}{t}=\mathrm{E}[N_{\bold{Y}}(0,1]],\quad\text{a.s..}
\]

\noindent By the above definition of Poisson cluster process it is
clear that $\X$ is ergodic with finite intensity $\nu\EE[S]$.~In
particular,

\begin{equation}\label{lln0}
\lim_{t\rightarrow\infty}\frac{N_{\bold{X}}(0,t]}{t}=\nu\mathrm{E}[S],\quad\text{a.s..}
\end{equation}

\subsection{Hawkes processes}\label{sec:prel1}
We say that $\bold{X}\subset\mathbb{R}$ is a Hawkes process if it
is a Poisson cluster process with $(\mathrm{b})$ in the definition
above replaced by:

\begin{itemize}

\item[(b)'] Each immigrant $X_i$ generates a cluster
$C_i=C_{X_i}$, which is the random set formed by the points of
generations $n=0,1,\ldots$ with the following branching structure:
the immigrant $X_i$ is said to be of generation 0.~Given
generations $0,1,\ldots,n$ in $C_i$, each point $Y\in C_i$ of
generation $n$ generates a Poisson process on $(Y,\infty)$, say
$\Phi$, of offspring of generation $n+1$ with intensity function
$h(\cdot-Y)$.~Here $h:(0,\infty)\rightarrow [0,\infty)$ is a
non-negative Borel function called fertility rate.

\end{itemize}

\noindent We refer the reader to Section 2 in M{\o}ller and
Rasmussen (2005) for more insight into the branching structure and
self-similarity property of clusters.~Consider the mean number of
points in any offspring process $\Phi$:

\[
\mu=\int_{0}^{\infty}h(t)\,\mathrm{d}t.
\]

\noindent As usual in the literature on Hawkes processes,
throughout this paper we assume

\begin{equation}\label{stability}
0<\mu<1.
\end{equation}

\noindent Condition $\mu>0$ excludes the trivial case in which
there are almost surely no offspring.~Recalling that the total
number of points in a cluster is equivalent to the total progeny
of the Galton-Watson process with one ancestor and number of
offspring per individual following a Poisson distribution with
mean $\mu$ (see p.~496 of Hawkes and Oakes (1974)), the other
condition $\mu<1$ is equivalent to assuming that
$\mathrm{E}[S]=1/(1-\mu)<\infty$.~For our purposes it is important
to recall that for Hawkes processes the distribution of $S$ is
given by

\begin{equation}\label{distrs}
P(S=k)=\frac{\mathrm{e}^{-k\mu}(k\mu)^{k-1}}{k!},\quad{\text{$k=1,2,\ldots$}}
\end{equation}

\noindent This follows by Theorem 2.11.2 in the book by Jagers
(1975).~Finally, since $\bold{X}$ is ergodic with a finite and
positive intensity equal to $\nu/(1-\mu)$ it holds:

\begin{equation}\label{lln}
\lim_{t\rightarrow\infty}\frac{N_{\bold{X}}(0,t]}{t}=\frac{\nu}{1-\mu},\quad\text{a.s..}
\end{equation}

\subsection{Large deviation principles}
\label{sec:prel2} We recall here some basic definitions in large
deviations theory (see, for instance, the book by Dembo and
Zeitouni (1998)).~A family of probability measures $\{\mu
_\alpha\}_{\alpha\in (0,\infty)}$ on a topological space
$(M,\mathcal{T}_M)$ satisfies the large deviations principle (LDP
for short) with rate function $J(\cdot)$ and speed $v(\cdot)$ if
$J:M\rightarrow [0,\infty ]$ is a lower semi-continuous function,
$v:[0,\infty)\rightarrow [0,\infty)$ is a measurable function
which increases to infinity, and the following inequalities hold
for every Borel set $B$:

\[
-\inf_{x\in B^\circ }J(x)\leq \liminf _{\alpha\rightarrow
\infty}\frac{1}{v(\alpha)}\log\mu_\alpha (B)\leq\limsup
_{\alpha\rightarrow\infty}\frac{1}{v(\alpha)}\log\mu_{\alpha}
(B)\leq -\inf _{x\in \overline{B}}J(x),
\]

\noindent where $B^\circ$ denotes the interior of $B$ and
$\overline{B}$ denotes the closure of $B$.~Similarly, we say that
a family of $M$-valued random variables $\{V_\alpha\}_{\alpha\in
(0,\infty)}$ satisfies the LDP if $\{\mu_\alpha\}_{\alpha\in
(0,\infty)}$ satisfies the LDP and $\mu_\alpha
(\cdot)=P(V_\alpha\in\cdot)$. We point out that the lower
semi-continuity of $J(\cdot)$ means that its level sets:

\[
\{x\in M: J(x)\leq a\},\quad\text{$a\geq 0$,}
\]

\noindent are closed; when the level sets are compact the rate
function $J(\cdot)$ is said to be good.

\section{Scalar large deviations}

\subsection{Scalar large deviations of Poisson cluster processes}\label{sec:LDPscalarPC}
Consider the ergodic Poisson cluster process $\X$ described
above.~In this section we prove that the process
$\{N_{\bold{X}}(0,t]/t\}$ satisfies a LDP on $\mathbb{R}$.~Define
the set

\[
\D_S=\{\theta\in\R:\mathrm{E}[\mathrm{e}^{\theta S} ]<\infty\}.
\]

\noindent With a little abuse of notation, denote by $C_0$ the
cluster generated by an immigrant at $0$ and let $L=\sup_{Y\in
C_0}|Y|$ be the radius of $C_0$.~We shall consider the following
conditions:

\begin{equation}\label{eq:hypPC}
\text{the function}\quad\theta\mapsto\mathrm{E}[\mathrm{e}^{\theta
S}]\quad\text{is essentially smooth and $0\in\D_S^{\circ}$}
\end{equation}

\noindent and

\begin{equation}\label{eq:radCo}
\mathrm{E}[L\mathrm{e}^{\theta S}]<\infty\quad\text{for all
$\theta\in\D_S^{\circ}$.}
\end{equation}

\noindent For the definition of essentially smooth function, we
refer the reader to Definition 2.3.5. in Dembo and Zeitouni
(1998).

\begin{Remark} \textnormal{
Since $S\geq 1$ we have that the function
$\varphi(\theta)=\mathrm{E}[\mathrm{e}^{\theta S}]$ is
increasing.~It follows that $\D_S^{\circ}=(-\infty,\theta_0)$ with
$\theta_0\in [0,\infty]$.~By the dominated convergence theorem we
have that $\varphi'(\theta)=\mathrm{E}[S\mathrm{e}^{\theta S}]$
and $\varphi''(\theta)=\mathrm{E}[S^2\mathrm{e}^{\theta S}]$, for
all $\theta\in\D_S^{\circ}$.~Hence, if $\theta_0<\infty$, to prove
that $\varphi$ is essentially smooth it suffices to show that
$\mathrm{E}[S\mathrm{e}^{\theta_0 S}]=\infty$.~On the other hand,
if $\theta_0=+\infty$, the function $\varphi$ is always
essentially smooth.}
\end{Remark}

\noindent It holds:

\begin{Theorem}\label{scalarLDPPC}
Assume $(\ref{eq:hypPC})$ and $(\ref{eq:radCo})$.~Then
$\{N_{\bold{X}}(0,t]/t\}$ satisfies a LDP on $\mathbb{R}$ with
speed $t$ and good rate function

\begin{equation}\label{rfPC}
\Lambda^*(x)=\sup_{\theta\in\R}(\theta x-\Lambda(\theta)).
\end{equation}

\noindent where $\Lambda (\theta)=\nu(\EE[\mathrm{e}^{\theta S}]-
1)$.
\end{Theorem}

\noindent It is easily verified that
$\Lambda^*(\nu\mathrm{E}[S])=0$.~Moreover, this is the unique zero
of $\Lambda^*(\cdot)$.~Therefore the probability law of
$N_{\bold{X}}(0,t]/t$ concentrates in arbitrarily small
neighborhoods of $\nu\mathrm{E}[S]$ as $t\rightarrow\infty$, as
stated by the law of large numbers $(\ref{lln0})$.~The LDP is a
refinement of the law of large numbers in that it gives us the
probability of fluctuations away the most probable value.

Before proving Theorem \ref{scalarLDPPC} we show that the same LDP
holds for the non-stationary Poisson cluster process ${\bf
X}_{t,T}$ with immigrant process empty on
$(-\infty,-T)\cup(t+T,\infty)$, where $T>0$ is a fixed
constant.~Furthermore, the LDP for ${\bf X}_{t,T}$ holds under a
weaker condition.

\begin{Theorem}\label{LDPfortrucatedPC}
Assume $(\ref{eq:hypPC})$.~Then $\{N_{{\bf X} _{t,T}}(0,t]/t\}$
satisfies a LDP on $\mathbb{R}$ with speed $t$ and good rate
function $(\ref{rfPC})$.
\end{Theorem}

\noindent$\bold{Proof}.$ The proof is based on the
G{\"a}rtner-Ellis theorem (see, for instance, Theorem 2.3.6 in
Dembo and Zeitouni (1998)).~We start proving that

\begin{equation}\label{le:gartner-ellis limitPC}
\lim_{t\rightarrow\infty}\frac{1}{t}\log\mathrm{E}[\mathrm{e}^{\theta
N_{{\bf X} _{t,T}}(0,t]}]=\left \{
\begin{array}{ll}
\nu(\mathrm{E}[\mathrm{e}^{\theta S}]-1)&\ {\rm if}\ \theta \in \D_S\\
+\infty &\ {\rm if}\ \theta\notin \D_S
\end{array}
\right .
\end{equation}

\noindent For a Borel set $A \subset \R$, let $I_{|A}=I\cap A$ be
the point process of immigrants in $A$.~Clearly $I_{|(0,t]}$,
$I_{|[-T,0]}$ and $I_{|(t,t+T]}$ are independent Poisson processes
with intensity $\nu$, respectively on $(0,t]$, $[-T,0]$ and
$(t,T+t]$.~Since $I_{|(0,t]}$, $I_{|[-T,0]}$ , and $I_{|(t,t+T]}$
are independent, by the definition of Poisson cluster process it
follows that the random sets $\{C_i:X_i\in I_{|(0,t]}\}$,
$\{C_i:X_i\in I_{|[-T,0]}\}$ and $\{C_i:X_i\in I_{|(t,t+T]}\}$ are
independent.~Therefore, for all $\theta\in\mathbb{R}$,

\begin{align}
\mathrm{E}\left[\mathrm{e}^{\theta N_{{\bf X}
_{t,T}}(0,t]}\right]&=
\mathrm{E}\left[\mathrm{e}^{\theta\left(\sum_{X_i\in
I_{|(0,t]}}N_{C_i}(0,t]+\sum_{X_i\in I_{|[-T,0]}} N_{C_i}(0,t] + \sum_{X_i\in I_{|(t,t+T]}}N_{C_i}(0,t]\right)}\right]
\nonumber\\
&=\mathrm{E}\left[\mathrm{e}^{\theta\sum_{X_i\in
I_{|(0,t]}}N_{C_i}(0,t]}\right]\mathrm{E}\left[\mathrm{e}^{\theta\sum_{X_i\in
I_{|[-T,0]}}N_{C_i}(0,t]}\right]\EE
\left[\mathrm{e}^{\theta\sum_{X_i\in
I_{|(t,t+T]}}N_{C_i}(0,t]}\right]. \nonumber
\end{align}

\noindent We shall show

\begin{equation}\label{limit2PC}
\lim_{t\rightarrow\infty}\frac{1}{t}\log\mathrm{E}\left[\mathrm{e}^{\theta\sum_{X_i\in
I_{|(0,t]}}N_{C_i}(0,t]}\right]=\left\{
\begin{array}{ll}
\nu(\mathrm{E}[\mathrm{e}^{\theta S}]-1)&\ {\rm if}\ \theta \in \D_S\\
+\infty &\ {\rm if}\ \theta\notin \D_S
\end{array}
\right .
\end{equation}

\noindent and

\begin{equation}\label{limit1PC}
\lim_{t\rightarrow\infty}\frac{1}{t}\log\mathrm{E}\left[\mathrm{e}^{\theta\sum_{X_i\in
I_{|[-T,0]}}N_{C_i}(0,t]}\right]=
\lim_{t\rightarrow\infty}\frac{1}{t}\log\mathrm{E}\left[\mathrm{e}^{\theta\sum_{X_i\in
I_{|(t,t+T]} }N_{C_i}(0,t]}\right] = 0,\quad\text{for $\theta\in
\D_S$.}
\end{equation}

\noindent Note that (\ref{le:gartner-ellis limitPC}) is a
consequence of $(\ref{limit2PC})$ and $(\ref{limit1PC})$. We first
prove $(\ref{limit2PC})$.~With a little abuse of notation, denote
by $C_0$ the cluster generated by an immigrant at $0$.~Since
$\{(X_i,C_i):X_i\in I_{|(0,t]}\}$ is an independently marked
Poisson process, by Lemma 6.4.VI in Daley and Vere-Jones (2003) we
have

\begin{align}
\mathrm{E}\left[\mathrm{e}^{\theta\sum_{X_i\in
I_{|(0,t]}}N_{C_i}(0,t]}\right]&=
\mathrm{E}\left[\mathrm{e}^{\theta\sum_{X_i\in I_{|(0,t]}}N_{C_i-X_i}(-X_i,t-X_i]}\right]\nonumber\\
&=\exp\left(\nu\int_{0}^{t}\mathrm{E}\left[\mathrm{e}^{\theta
N_{C_0}(-x,t-x]}-1\right]\,\mathrm{d}x\right) \nonumber
\\
&=\exp\left(\nu t \int_{0}^{1}\mathrm{E}\left[\mathrm{e}^{\theta
N_{C_0}(-tz ,(1-z)
t]}-1\right]\,\mathrm{d}z\right).\label{lintsolPC}
\end{align}

\noindent Therefore if $\theta\in\D_S$, the expectation in
(\ref{lintsolPC}) goes to $\mathrm{E}[\mathrm{e}^{\theta S}-1]$ as
$t\rightarrow\infty$ by the monotone convergence theorem.~Hence,
for $\theta \in \D_S$ the limit (\ref{limit2PC}) follows from the
dominated convergence theorem.~For $\theta\notin\D_S$ the
expectation in $(\ref{lintsolPC})$ goes to $+\infty$ as
$t\rightarrow\infty$ by the monotone convergence theorem, and the
limit $(\ref{limit2PC})$ follows by Fatou's lemma.~We now show
$(\ref{limit1PC})$.~Here again, since $\{(X_i,C_i):X_i\in
I_{|[-T,0]}\}$ is an independently marked Poisson process, by
Lemma 6.4.VI in Daley and Vere-Jones (2003) we have

\begin{align}
\mathrm{E}\left[\mathrm{e}^{\theta\sum_{X_i\in
I_{|[-T,0]}}N_{C_i}(0,t]}\right]&=\mathrm{E}\left[\mathrm{e}^{\theta\sum_{X_i\in
I_{|[-T,0]}}N_{C_i-X_i}(-X_i,t-X_i]}\right]\nonumber\\
&=\exp\left(\nu\int_{0}^{T}\mathrm{E}\left[\mathrm{e}^{\theta
N_{C_0}(x,x+t]}-1\right]\,\mathrm{d}x\right).
\end{align}

\noindent Now note that, for $\theta \in \D_S \cap [0,\infty)$, we
have

\[
0\leq\frac{1}{t}\log\mathrm{E}\left[\mathrm{e}^{\theta\sum_{X_i\in
I_{|[-T,0]}}N_{C_i}(0,t]}\right]\leq\frac{\nu}{t}\int_{0}^{T}\mathrm{E}\left[\mathrm{e}^{\theta
S}-1\right]\,\mathrm{d}x <\infty
\]

\noindent and, for each $\theta\leq 0$,

\[
\frac{\nu T}{t}\mathrm{E}\left[\mathrm{e}^{\theta
S}-1\right]\leq\frac{\nu}{t}\int_{0}^{T}\mathrm{E}\left[\mathrm{e}^{\theta
N_{C_0}(x,x+t]}-1\right]\,\mathrm{d}x\leq 0.
\]

\noindent  By passing to the limit as $t\rightarrow\infty$ we get
that the first limit in $(\ref{limit1PC})$ is equal to 0.~The
proof for the second limit in $(\ref{limit1PC})$ is rigorously the
same.~Hence we proved (\ref{le:gartner-ellis limitPC}).~Using
assumption (\ref{eq:hypPC}), the conclusion is a consequence of
the G{\"a}rtner-Ellis theorem.
\\
$\bold{\square}$

\noindent $\bold{Proof\,of\,Theorem \,\ref{scalarLDPPC}}.$ The
proof is similar to that one of Theorem \ref{LDPfortrucatedPC} and
is again based on the G{\"a}rtner-Ellis theorem.~We start showing
that

\begin{equation}\label{eq:GE}
\lim_{t\rightarrow\infty}\frac{1}{t}\log\mathrm{E}[\mathrm{e}^{\theta
N_{{\bf X}}(0,t]}]=\left \{
\begin{array}{ll}
\nu(\mathrm{E}[\mathrm{e}^{\theta S}]-1)&\ {\rm if}\ \theta \in \D_S^{\circ}\\
+\infty &\ {\rm if}\ \theta\notin \D_S
\end{array}
\right .
\end{equation}

\noindent By similar arguments as in the proof of Theorem
\ref{LDPfortrucatedPC}, using the definition of ${\bf X}$, we have

\[
\mathrm{E}\left[\mathrm{e}^{\theta N_{{\bf
X}}(0,t]}\right]=\mathrm{E}\left[\mathrm{e}^{\theta N_{{\bf X}
_{t,T}}(0,t]}\right]\mathrm{E}\left[\mathrm{e}^{\theta (N_{{\bf
X}}(0,t]-N_{{\bf X} _{t,T}}(0,t])}\right],\quad\text{for all
$\theta\in\R$, $t>0$.}
\]

\noindent By the computations in the proof of Theorem
\ref{LDPfortrucatedPC}, in order to prove (\ref{eq:GE}) we only
need to check that

\begin{equation}\label{eq:finallimit}
\lim_{t\to
\infty}\frac{1}{t}\log\mathrm{E}\left[\mathrm{e}^{\theta (N_{{\bf
X}}(0,t]-N_{{\bf X} _{t,T}}(0,t])}\right]=0,\quad\text{for all
$\theta\in\D_S^\circ$.}
\end{equation}

\noindent It is easily verified for $\theta\leq 0$ (the
argument of the expectation is bounded below by
$\mathrm{e}^{\theta S}$ and above by $1$). We only check
$(\ref{eq:finallimit})$ for $\theta\in\D_S^{\circ}\cap
(0,\infty)$.~Here again, for a Borel set $A\subset\R$, let
$I_{|A}=I\cap A$ denote the point process of immigrants in
$A$.~Note that

\[
N_{{\bf X}}(0,t]-N_{{\bf X} _{t,T}}(0,t]=\sum_{X_i\in
I_{|(-\infty,-T)}}N_{C_i}(0,t]+\sum_{X_i\in I_{|(t+T,\infty)}}
N_{C_i}(0,t],\quad\text{$t>0$.}
\]

\noindent Clearly $I_{|(-\infty,-T)}$ and $I_{|(t+T,\infty)}$ are
independent Poisson processes with intensity $\nu$, respectively
on $(-\infty,-T)$ and $(t+T,\infty)$.~Thus, by the definition of
Poisson cluster process it follows that the random sets
$\{C_i:X_i\in I_{|(-\infty,-T)}\}$ and $\{C_i:X_i\in
I_{|(t+T,\infty)}\}$ are independent.~Therefore, for all
$\theta\in\D_S^{\circ}\cap (0,\infty)$,

\begin{equation*}
\mathrm{E}\left[\mathrm{e}^{\theta (N_{{\bf X}}(0,t]-N_{{\bf X}
_{t,T}}(0,t])}\right]=\mathrm{E}\left[\mathrm{e}^{\theta\sum_{X_i\in
I_{|(-\infty,-T)}}N_{C_i}(0,t]}\right]\mathrm{E}\left[\mathrm{e}^{\theta\sum_{X_i\in
I_{|(t+T,\infty)}}N_{C_i}(0,t]}\right].
\end{equation*}

\noindent Since $\{(X_i,C_i):X_i\in I_{|(-\infty,-T)}\}$ and
$\{(X_i,C_i):X_i\in I_{|(t+T,\infty)}\}$ are independently marked
Poisson processes, by Lemma 6.4.VI in Daley and Vere-Jones (2003)
we have

\begin{align}
\mathrm{E}\left[\mathrm{e}^{\theta\sum_{X_i\in
I_{|(-\infty,-T)}}N_{C_i}(0,t]}\right]&=
\mathrm{E}\left[\mathrm{e}^{\theta\sum_{X_i\in I_{|(-\infty,-T)}}N_{C_i-X_i}(-X_i,t-X_i]}\right]\nonumber\\
&=\exp\left(\nu\int_{T}^{\infty}\mathrm{E}\left[\mathrm{e}^{\theta
N_{C_0}(x,t+x]}-1\right]\,\mathrm{d}x\right)\nonumber
\end{align}

\noindent and

\begin{align}
\mathrm{E}\left[\mathrm{e}^{\theta\sum_{X_i\in
I_{|(t+T,\infty)}}N_{C_i}(0,t]}\right]&=\exp\left(\nu\int_{t+T}^{\infty}\mathrm{E}\left[\mathrm{e}^{\theta
N_{C_0}(-x,t-x]}-1\right]\,\mathrm{d}x\right)\nonumber\\
&=\exp\left(\nu\int_{T}^{\infty}\mathrm{E}\left[\mathrm{e}^{\theta
N_{C_0}(-t-z,-z]}-1\right]\,\mathrm{d}z\right).\nonumber
\end{align}

\noindent Now notice that since $\theta>0$ we have
\[
\mathrm{e}^{\theta N_{C_0}(x,x+t]}-1\leq(\mathrm{e}^{\theta
N_{C_0}(\R)}-1)\bold{1}\{x\leq L\},\quad\text{for all $x\geq T$}
\]
\noindent and
\[
\mathrm{e}^{\theta N_{C_0}(-t-z,-z]}-1\leq(\mathrm{e}^{\theta
N_{C_0}(\R)}-1)\bold{1}\{z\leq L\},\quad\text{for all $z\geq T$.}
\]
\noindent Relation $(\ref{eq:finallimit})$ follows by assumption
$(\ref{eq:radCo})$ noticing that the above relations yield
\[
\mathrm{E}\left[\mathrm{e}^{\theta\sum_{X_i\in
I_{|(-\infty,-T)}}N_{C_i}(0,t]}\right]\leq\exp(\nu\mathrm{E}[L(\mathrm{e}^{\theta
S}-1)]),\quad\text{for all $\theta\in\D_S^{\circ}\cap (0,\infty)$, $t>0$.}
\]
\noindent and
\[
\mathrm{E}\left[\mathrm{e}^{\theta\sum_{X_i\in
I_{|(t+T,\infty)}}N_{C_i}(0,t]}\right]\leq\exp(\nu\mathrm{E}[L(\mathrm{e}^{\theta
S}-1)]),\quad\text{for all $\theta\in\D_S^{\circ}\cap (0,\infty)$, $t>0$.}
\]
Therefore, (\ref{eq:GE}) is proved.~Now, if $\D_S=\D_S^{\circ}$
then the claim is a consequence of the G{\"a}rtner-Ellis theorem
and assumption (\ref{eq:hypPC}).~It remains to deal with the case
$\D_S\neq\D_S^{\circ}$.~We shall show the large deviations upper
and lower bounds proving that for any sequence $\{t_n\}_{n\geq
1}\subset (0,\infty)$ diverging to $+\infty$, as $n\to\infty$,
there exists a subsequence $\{s_n\}\subseteq\{t_n\}$ such that

\begin{equation}\label{eq:LDUB}
\limsup_{n\to\infty}\frac{1}{s_n}\log P(N_{\bf X}(0,s_n]/s_{n}\in
F)\leq -\inf_{x\in F}\Lambda^{*}(x),\quad\text{for all closed sets
$F$}
\end{equation}

\noindent and

\begin{equation}\label{eq:LDLB}
\liminf_{n\to\infty}\frac{1}{s_n}\log P(N_{\bf X}(0,s_n]/s_n\in
G)\geq -\inf_{x\in G}\Lambda^{*}(x),\quad\text{for all open sets
$G$,}
\end{equation}

\noindent where $\Lambda^*$ is defined by $(\ref{rfPC})$ (then the
large deviations upper and lower bounds hold for any sequence
$\{t_n\}$ and the claim follows).~By assumption (\ref{eq:hypPC}),
there exists $\theta_0>0$ such that $\D_S =
(\infty,\theta_0]$.~Let $\{t_n\}_{n\geq 1}\subset (0,\infty)$ be a
sequence diverging to $+\infty$, as $n\to\infty$, and define the
extended non-negative real number $l\in [0,\infty]$ by

\[
l\equiv\limsup_{n\to\infty}\frac{1}{t_n}\log\mathrm{E}[\mathrm{e}^{\theta_0
N_{\bf X}(0,t_n]}].
\]

\noindent Clearly, there exists a subsequence
$\{s_n\}\subseteq\{t_n\}$ which realizes this $\limsup$, i.e.

\[
\lim_{n\to\infty}\frac{1}{s_n}\log\mathrm{E}[\mathrm{e}^{\theta_0
N_{\bf X}(0,s_n]}]=l.
\]

\noindent By $(\ref{eq:GE})$ it follows that

\[
\lim_{n\to\infty}\frac{1}{s_n}\log\mathrm{E}[\mathrm{e}^{\theta
N_{\bf X}
(0,s_n]}]=\tilde{\Lambda}(\theta),\quad\text{$\theta\in\R$}
\]

\noindent where

\begin{equation*}
\tilde{\Lambda}(\theta)=\left \{
\begin{array}{ll}
\Lambda(\theta) &\ {\rm if}\ \theta<\theta_0\\
l &\ {\rm if}\ \theta=\theta_0\\
+\infty &\ {\rm if}\ \theta>\theta_0.
\end{array}
\right.
\end{equation*}

\noindent Note that, irrespective to the value of $l$,
$\tilde{\Lambda}$ is essentially smooth (however it may be not
lower semi-continuous).~We now show that the Legendre transform of
$\Lambda$ and $\tilde{\Lambda}$ coincide, i.e.

\begin{equation}\label{eq:tLtrasf}
\tilde{\Lambda}^{*}(x)=\Lambda^{*}(x),\quad\text{$x\in\R$}.
\end{equation}

\noindent A straightforward computation gives
$\tilde{\Lambda}^*(x)=\Lambda^*(x)= +\infty$, for $x<0$, and
$\tilde{\Lambda}^*(0)=\Lambda^*(0)=\nu$.~Now, note that since
$\theta_0<\infty$, $\tilde{\Lambda}^*(x)$ and $ \Lambda^*(x)$ are
both finite, for $x>0$.~Moreover, since $\Lambda$ and
$\tilde{\Lambda}$ are essentially smooth, if $x>0$ we have that
$\Lambda^*(x)=\theta_{x}x-\Lambda(\theta_x)$ and
$\tilde{\Lambda}^*(x)=\tilde \theta_{x}x-\tilde{\Lambda}(\tilde
\theta_x)$, where $\theta_x$ (respectively $\tilde \theta_x$) is
the unique solution of $\Lambda'(\theta) = x $ (respectively
$\tilde \Lambda ( \theta) = x$)  on $(-\infty,\theta_0)$.~The
claim $(\ref{eq:tLtrasf})$ follows recalling that
$\tilde{\Lambda}(\theta)=\Lambda(\theta)=\nu
(\mathrm{E}[\mathrm{e}^{\theta S}]-1)$ on $\D_S^{\circ}$.~Now,
applying part $(a)$ of Theorem 2.3.6 in Dembo and Zeitouni (1998)
we have $(\ref{eq:LDUB})$.~Applying part $(b)$ of Theorem 2.3.6 in
Dembo and Zeitouni (1998) we get

\begin{equation}\label{eq:tGEa}
\liminf_{n\to\infty}\frac{1}{s_n}\log P(N_{\bf X}(0,s_n]/s_{n}\in
G)\geq -\inf_{x\in G\cap\mathcal{F}}\Lambda^{*}(x),\quad\text{for
any open set $G$,}
\end{equation}

\noindent where $\mathcal{F}$ is the set of exposed points of
$\Lambda^*$ whose exposing hyperplane belongs to
$(-\infty,\theta_0)$, i.e.

\[
\mathcal{F}=\{y\in\R:\exists\,\theta\in\D_S^{\circ}\hbox{ such
that for all } x \ne y, \theta y - \Lambda^* (y) > \theta x -
\Lambda^* (x)\}.
\]

\noindent We now prove that $\mathcal{F} = (0,+\infty)$.~For $y<0$,
$\Lambda^*(y)=\infty$, therefore an exposing hyperplane satisfying
the corresponding inequality does not exist.~For $y>0$ consider
the exposing hyperplane $\theta= \theta_y$, where $\theta_y$ is
the unique positive solution on $(-\infty,\theta_0)$ of
$\mathrm{E}[S\mathrm{e}^{\theta S}]=y/\nu$.~Note that
$\Lambda'(\theta)=\mathrm{E}[S\mathrm{e}^{\theta S}]$ and
$\Lambda''(\theta)=\mathrm{E}[S^2\mathrm{e}^{\theta S}]$ for all
$\theta<\theta_0$.~In particular, since $S\geq 1$, we have that
$\Lambda$ is strictly convex on $(-\infty,\theta_0)$.~Therefore,
for all $x\ne y$, it follows
\begin{align}
\theta_{y}y-\Lambda^*(y)=\Lambda(\theta_y)&>\Lambda(\theta_x)+\Lambda'(\theta_x)(\theta_y-\theta_x)\nonumber\\
&=\theta_{y}x-\Lambda^{*}(x).\nonumber
\end{align}
\noindent It remains to check that $0 \notin \mathcal{F}$.~Notice
that since $\mathrm{E}[S\mathrm{e}^{\theta_x S}] = x / \nu$,
$\lim_{x\downarrow 0} \theta_x=-\infty$.~Also, by the implicit
function theorem, $ x \mapsto \theta_x$ is a continuous mapping on
$(0,\infty)$.~Now assume that $0 \in \mathcal{F}$, then there
would exist $\theta < \theta_0$, such that for all $x >0$, $-
\Lambda^*(0) > \theta x - \Lambda^*(x)$.~However, by the
intermediate values theorem, there exists $y>0$ such that $\theta
=\theta_y$, and we obtain a contradiction.~This implies
$\mathcal{F} = (0,+\infty) $ as claimed.~Now recall that
$\Lambda^*(x)=+\infty$ for $x<0$; moreover, $\lim_{x\downarrow
0}\Lambda^{*}(x)=\Lambda^*(0)=\nu$ (indeed, $\lim_{x\downarrow
0}\theta_x=-\infty$).~Therefore
\[
\inf_{x\in G\cap\mathcal{F}}\Lambda^{*}(x)\leq\inf_{x\in
G}\Lambda^{*}(x),\quad\text{for any open set $G$.}
\]
\noindent Finally, by $(\ref{eq:tGEa})$ and the above inequality
we obtain $(\ref{eq:LDLB})$.\\
\noindent$\bold{\square}$

\subsection{Scalar large deviations of Hawkes processes}
\label{sec:LDPscalar} Consider the ergodic Hawkes process $\X$
described before.~In this section we prove that the process
$\{N_{\bold{X}}(0,t]/t\}$ satisfies a LDP, and we give the
explicit expression of the rate function.~Our result is a
refinement of the law of large numbers $(\ref{lln})$.~The
following theorem holds:

\begin{Theorem}\label{scalarLDP}
Assume $(\ref{stability})$ and

\begin{equation}\label{eq:hfinmean}
\int_0^\infty th(t)\,\mathrm{d}t<\infty.
\end{equation}

\noindent Then $\{N_{\bold{X}}(0,t]/t\}$ satisfies a LDP on
$\mathbb{R}$ with speed $t$ and good rate function

\begin{equation}\label{rf}
\Lambda^*(x)=\left \{
\begin{array}{ll}
x\theta_x+\nu-\frac{\nu x}{\nu+\mu x}&\ {\rm if}\ x\in (0,\infty) \\
\nu &\ {\rm if}\ x=0\\
+\infty &\ {\rm if}\ x\in (-\infty,0)
\end{array}
\right. ,
\end{equation}

\noindent where $\theta=\theta_x$ is the unique solution in
$(-\infty,\mu-1-\log\mu)$ of

\begin{equation}\label{eq:galtonprogeny}
\mathrm{E}\left[S\mathrm{e}^{\theta
S}\right]=x/\nu,\quad\text{$x>0$,}
\end{equation}

\noindent or equivalently of

\[
\mathrm{E}[\mathrm{e}^{\theta
S}]=\frac{x}{\nu+x\mu},\quad\text{$x>0$.}
\]
\end{Theorem}

\vspace{0.5cm}

\noindent$\bold{Proof}.$ The proof is a consequence of Theorem
\ref{scalarLDPPC}.~We start noticing that by $(\ref{distrs})$ we
have

\[
\mathrm{E}[\mathrm{e}^{\theta S}]=\sum_{k\geq 1}\frac{(
\mathrm{e}^{\theta-\mu})^k(k\mu)^{k-1}}{k!},
\]

\noindent and this sum is infinity for $\theta>\mu-1-\log\mu$ and
finite for $\theta<\mu-1-\log\mu$ (apply, for instance, the ratio
criterion).~If $\theta=\mu-1-\log\mu$ the sum above is
finite.~Indeed, in this case

\[
\mathrm{E}[\mathrm{e}^{\theta S}]=(1/\mu)\sum_{k\geq 1}\frac{
\mathrm{e}^{-k}k^{k-1}}{k!}=1/\mu.
\]

\noindent Therefore $\D_{S}=(-\infty,\mu-1-\log\mu]$.~The origin
belongs to $\mathcal{D}_{S}^{\circ}$ in that by
$(\ref{stability})$ and the inequality $\mathrm{e}^{x}>x+1$,
$x\neq 0$, we have $\frac{\mathrm{e}^{\mu-1}}{\mu}>1$.~The
function $\theta\mapsto \EE[\mathrm{e}^{\theta S}]$ is essentially
smooth.~Indeed, it is differentiable in the interior of $\D_S$ and

\[
\mathrm{E}[S\mathrm{e}^{(\mu-1-\log\mu)S}]=\infty
\]

\noindent because

\[
\mathrm{E}[S\mathrm{e}^{(\mu-1-\log\mu)S}]=(1/\mu)\sum_{k\geq
1}\frac{ \mathrm{e}^{-k}k^{k}}{k!}
\]

\noindent and this sum is infinity since by Stirling's formula
$\frac{\mathrm{e}^{-k}k^{k}}{k!}\sim 1/\sqrt{2\pi k}$.~We now
check assumption $(\ref{eq:radCo})$.~By the structure of the
clusters, it follows that there exists a sequence of independent
non-negative random variables $\{V_n\}_{n\geq 1}$, independent of
$S$, such that $V_1$ has probability density $h(\cdot)/\mu$ and
the following stochastic domination holds:

\[
L\leq\sum_{n=1}^{S}V_n,\quad\text{a.s.}
\]

\noindent (see Reynaud-Bouret and Roy (2007)).~Therefore, for all
$\theta<\mu-1-\log\mu$, we have

\[
\mathrm{E}[L\mathrm{e}^{\theta
S}]\leq\mathrm{E}\left[\mathrm{e}^{\theta
S}\sum_{n=1}^{S}V_n\right]=\mathrm{E}[V_1]\mathrm{E}[S\mathrm{e}^{\theta
S}].
\]

\noindent Since $\theta<\mu-1-\log \mu$, we have
$\mathrm{E}[S\mathrm{e}^{\theta S}] < \infty$; moreover,
assumption $(\ref{eq:hfinmean})$ yields
$$
\mathrm{E}[V_1]=\frac{1}{\mu}\int_{0}^{\infty}th(t)\,\mathrm{d}t<
\infty.
$$
Hence, condition $(\ref{eq:radCo})$ holds, and by Theorem
\ref{scalarLDPPC}, $\{N_{\bold{X}}(0,t]/t\}$ satisfies a LDP on
$\mathbb{R}$ with speed $t$ and good rate function

\[
\Lambda^*(x)=\sup_{\theta\in\mathbb{R}}(\theta
x-\Lambda(\theta))=\sup_{\theta\leq\mu-1-\log\mu}(\theta
x-\Lambda(\theta)).
\]

\noindent Now $\Lambda^*(x)=\infty$ if $x<0$, in that in such a
case $\lim_{\theta\rightarrow-\infty}(\theta
x-\Lambda(\theta))=\infty$.~If $x>0$, letting $\theta_x\in
(-\infty,\mu-1-\log\mu)$ denote the unique solution of the
equation $(\ref{eq:galtonprogeny})$ easily follows that

\begin{equation}\label{eq:legendre}
\Lambda^*(x)=x\theta_x-\Lambda(\theta_x).
\end{equation}

\noindent It is well-known (see, for instance, p. 39 in Jagers
(1975)) that, for all $\theta\in (-\infty,\mu-1-\log\mu)$,
$\mathrm{E}[\mathrm{e}^{\theta S}]$ satisfies

\[
\mathrm{E}[\mathrm{e}^{\theta
S}]=\mathrm{e}^{\theta}\exp\{\mu(\mathrm{E}[\mathrm{e}^{\theta
S}]-1)\},
\]

\noindent therefore differentiating with respect to $\theta$ we
get

\begin{equation}\label{pasqua}
\mathrm{E}[S\mathrm{e}^{\theta
S}]=\frac{\mathrm{e}^{\theta}\exp\{\mu(\mathrm{E}[\mathrm{e}^{\theta
S}]-1)\}}{1-\mu\mathrm{e}^{\theta}\exp\{\mu(\mathrm{E}[\mathrm{e}^{\theta
S}]-1)\}}=\frac{\mathrm{E}[\mathrm{e}^{\theta
S}]}{1-\mu\mathrm{E}[\mathrm{e}^{\theta S}]}.
\end{equation}

\noindent Setting $\theta=\theta_x$ in the above equality and
using $(\ref{eq:galtonprogeny})$ we have

\[
\frac{x}{\nu}=\frac{\mathrm{E}[\mathrm{e}^{\theta_x
S}]}{1-\mu\mathrm{E}[\mathrm{e}^{\theta_x S}]},
\]

\noindent which yields

\[
\mathrm{E}[\mathrm{e}^{\theta_x S}]=\frac{x}{\nu+x\mu}.
\]

\noindent Thus, by $(\ref{eq:legendre})$ we have for $x>0$

\[
\Lambda^*(x)=x\theta_x+\nu-\frac{\nu x}{\nu+\mu x}.
\]

\noindent The conclusion follows noticing that a direct computation gives $\Lambda^*(0)=\nu$.\\
\noindent$\bold{\square}$

\section{Sample path large deviations}\label{sec:SP-LDP}
Let $\X$ be the ergodic Poisson cluster process described at the
beginning.~The results proved in this section are sample path LDP
for $\X$.

\subsection{Sample path large deviations in the topology of point-wise convergence}
Let $D[0,1]$ be the space of c\`adl\`ag functions on the interval
$[0,1]$.~Here we prove that
$\{\frac{N_{\bold{X}}(0,\alpha\cdot]}{\alpha}\}$ satisfies a LDP
on $D[0,1]$ equipped with the topology of point-wise convergence
on $D[0,1]$.~The LDP we give
is a refinement of the following functional law of
large numbers:

\begin{equation}\label{flln}
\lim_{\alpha\rightarrow\infty}\frac{N_{\bold{X}}(0,\alpha
\cdot]}{\alpha}=\chi(\cdot)\quad\text{a.s.,}
\end{equation}

\noindent where $\chi(t)=\nu\mathrm{E}[S]t$.~As this is a
corollary of the LDP we establish, we do not include a separate
proof of this result.~Letting $\Lambda^*(\cdot)$ denote the rate
function of the scalar LDP, we have:

\begin{Theorem}\label{SPLDP}
Assume $(\ref{eq:hypPC})$ and $(\ref{eq:radCo})$.~If moreover
$\mathcal{D}_S$ is open, then $\{\frac{N_{\bold{X}}(0,\alpha
\cdot]}{\alpha}\}$ satisfies a LDP on $D[0,1]$, equipped with the
topology of point-wise convergence, with speed $\alpha$ and good
rate function

\begin{equation}
\label{rfIexpression} J(f)=\left \{
\begin{array}{ll}
\int _0^1 \Lambda ^*(\dot{f}(t))dt&\ {\rm if} f\in AC_0[0,1] \\
\infty &\ otherwise \\
\end{array}
\right. ,
\end{equation}

\noindent where $AC_0[0,1]$ is the family of absolutely continuous
functions $f(\cdot)$ defined on $[0,1]$, with $f(0)=0$.
\end{Theorem}

\noindent While it is tempting to conjecture that the result above
holds even if the effective domain of $S$ is not open, we do not
have a proof of this claim.~If we take
$\chi(t)=\nu\mathrm{E}[S]t$, then $J(\chi)=0$.~Moreover this is
the unique zero of $J(\cdot)$.~Thus the law of
$N_{\bold{X}}(0,\alpha\cdot]/\alpha$ concentrates in arbitrarily
small neighborhoods of $\chi(\cdot)$ as $\alpha\rightarrow\infty$,
as ensured by the functional law of large numbers
$(\ref{flln})$.~The sample path LDP is a refinement of the
functional law of large numbers in that it gives the probability
of fluctuations away the most likely path.

As in Section \ref{sec:LDPscalarPC}, denote by ${\bf X}_{t,T}$ the
non-stationary Poisson cluster process with immigrant process
empty on $(-\infty,-T)\cup (t+T,\infty)$, where $T>0$ is a fixed
constant.~Before proving Theorem \ref{SPLDP} we show that the same
LDP holds for ${\bf X}_{t,T}$.~Furthermore, the LDP for ${\bf
X}_{t,T}$ holds under a weaker condition.

\begin{Theorem}\label{LDPfortrucated2}
Assume $(\ref{eq:hypPC})$.~Then
$\{N_{\bold{X}_{\alpha\cdot,T}}(0,\alpha\cdot]/\alpha\}$ satisfies
a LDP on $D[0,1]$, equipped with the topology of point-wise
convergence, with speed $\alpha$ and good rate function
$(\ref{rfIexpression})$.
\end{Theorem}

\noindent To prove this theorem we need Lemma \ref{le:agcmgt}
below, whose proof can be found in Ganesh, Macci and Torrisi
(2005) (see Lemma 2.3 therein).

\begin{Lemma}\label{le:agcmgt}
Let $(\theta_1,\ldots,\theta_n)\in\mathbb{R}^n$ and let
$w_1,\ldots,w_n\geq 0$ be such that $w_1\leq\ldots\leq w_n$.~Then
$\sum_{i=k}^{n}\theta_i w_i\leq\theta^*w^*$ for all
$k\in\{1,\ldots,n\}$, for any
$\theta^*\geq\max\{\max\{\sum_{i=k}^{n}\theta_i:k\in\{1,\ldots,n\}\},0\}$
and any $w^*\geq w_n$.
\end{Lemma}

\noindent$\bold{Proof\,of\,Theorem\,\ref{LDPfortrucated2}}.$ With
a little abuse of notation denote by $C_0$ the cluster generated
by an immigrant at $0$.~We first show the theorem under the
additional condition

\begin{equation}\label{eq:addcond}
N_{C_0}((-\infty,0))=0,\quad\text{a.s..}
\end{equation}

\noindent The idea in proving Theorem \ref{LDPfortrucated2} is to
apply the Dawson-G\"artner theorem to "lift" a LDP for the
finite-dimensional distributions of $\{N_{\bold{X}_{\alpha
t,T}}(0,\alpha t]/\alpha\}$ to a LDP for the process.~Therefore,
we first show the following claim:
\begin{itemize}
\item[($\bold{C}$)] For all $n\geq 1$ and $0\leq t_1<\ldots
<t_n\leq 1$, $\left(N_{\bold{X}_{\alpha t_1,T}}(0,\alpha
t_1]/\alpha, \ldots,N_{\bold{X}_{\alpha t_n,T}}(0,\alpha
t_n]/\alpha\right)$ satisfies the LDP in $\mathbb{R}^n$ with speed
$\alpha$ and good rate function
\begin{equation}
\label{It1-tnExpression} J_{t_1,\ldots ,t_n}(x_1,\ldots ,x_n)=
\sum_{j=1}^{n}(t_j-t_{j-1})\Lambda^{*}\left(\frac{x_j-x_{j-1}}{t_j-t_{j-1}}\right),
\end{equation}
where $x_0=0$ and $t_0=0$.
\end{itemize}

\noindent Claim $(\bold{C})$ is a consequence of the
G\"artner-Ellis theorem in $\mathbb{R}^n$, and will be shown in
three steps:
\begin{itemize}

\item[(a)] For each $(\theta_1,\ldots,\theta_n)\in\mathbb{R}^n$,
we prove that
\begin{equation}\label{limit(a)}
\Lambda_{t_1,\ldots ,t_n}(\theta_1,\ldots,\theta
_n)\equiv\lim_{\alpha\rightarrow\infty}\frac{1}{\alpha}\log\mathrm{E}
\left[\exp\left(\sum _{i=1}^n\theta_i N_{\bold{X}_{\alpha
t_i,T}}(0,\alpha t_i]\right)
\right]=\sum_{j=1}^{n}(t_j-t_{j-1})\Lambda\left(\sum
_{i=j}^n\theta_i\right),
\end{equation}
where the existence of the limit (as an extended real number) is
part of the claim, and $\Lambda(\cdot)$ is defined in the
statement of Theorem \ref{scalarLDPPC}.

\item[(b)] The function $\Lambda _{t_1,\ldots,t_n}(\cdot)$
satisfies the hypotheses of the G\"artner-Ellis theorem.

\item[(c)] The rate function

\[
J_{t_1,\ldots,t_n}(x_1,\ldots,x_n)\equiv\sup _{(\theta _1,\ldots
,\theta _n)\in\mathbb{R}^n}\left[\sum_{i=1}^{n}\theta_i x_i-
\Lambda_{t_1,\ldots,t_n}(\theta_1,\ldots,\theta_n)\right]
\]

coincides with the rate function defined in
$(\ref{It1-tnExpression})$.
\end{itemize}

\noindent{\em Proof of $(\mathrm a)$}.~For a Borel set
$A\subset\R$, denote by $I_{|A}=I\cap A$ the Poisson process of
immigrants in $A$.~Since, for each $t$, $I_{|(0,t]}$ and
$I_{|[-T,0]}$  are independent, it follows from the definition of
Poisson cluster process that, for each $i$, the random sets
$\{C_k:X_k\in I_{|(0,\alpha t_i]}\}$ and $\{C_k:X_k\in
I_{|[-T,0]}\}$ are independent.~Therefore,

\begin{align}
\label{eq:nlaplace} \mathrm{E}\left[\exp\sum_{i=1}^{n}\theta_i
N_{\bold{X}_{\alpha t_i,T}}(0,\alpha t_i]\right]
&=\mathrm{E}\left[\exp\sum_{i=1}^{n}\theta_i \sum_{X_k\in
I_{|(0,\alpha t_i]}}N_{C_k}(0,\alpha t_i]
\right]\mathrm{E}\left[\exp\sum_{i=1}^{n}\theta_i \sum_{X_k\in
I_{|[-T,0]}} N_{C_k}(0,\alpha t_i] \right],
\end{align}

\noindent where we used the independence and the assumption that
$N_{C_0}(-\infty,0)=0$ a.s..~In order to prove (\ref{limit(a)}),
we treat successively the two terms in
(\ref{eq:nlaplace}).~Viewing $I_{|(0,\alpha t_i]}$ as the
superposition of the $i$ independent Poisson processes:
$I_{|(\alpha t_{j-1},\alpha t_j]}$ on $(\alpha t_{j-1},\alpha
t_{j}]$ $(j=1,\ldots,i)$ with intensity $\nu$ we get

\begin{align}
\mathrm{E}\left[\exp\sum_{i=1}^{n}\theta_i \sum_{X_k\in
I_{|(0,\alpha t_i]}}N_{C_k}(0,\alpha t_i]
\right]&=\mathrm{E}\left[\exp\sum_{i=1}^{n}\sum_{j=1}^{i}
\sum_{X_k\in I_{|(\alpha t_{j-1},\alpha t_j]}}\theta_i
N_{C_k}(0,\alpha t_i]\right]\nonumber\\
&=\mathrm{E}\left[\exp\sum_{j=1}^{n}\sum_{i=j}^{n} \sum_{X_k\in
I_{|(\alpha t_{j-1},\alpha t_j]}}\theta_i
N_{C_k}(0,\alpha t_i] \right]\nonumber\\
&=\prod_{j=1}^{n}\mathrm{E}\left[\exp\sum_{i=j}^{n} \sum_{X_k\in
I_{|(\alpha t_{j-1},\alpha t_j]}}\theta_i N_{C_k}(0,\alpha
t_i]\right],\label{prodotto}
\end{align}

\noindent where in the latter equality we used the independence of
$\{C_k:X_k\in I_{|(\alpha t_{j-1},\alpha t_j]}\}$
$(j=1,\ldots,n)$.~Since, for each $j$, $\{(X_k,C_k):\,X_k\in
I_{|(\alpha t_{j-1},\alpha t_j]}\}$ is an independently marked
Poisson process, by Lemma 6.4.VI in Daley and Vere-Jones (2003) we
have

\begin{align}
&\mathrm{E}\left[\exp\sum_{i=j}^{n}\sum_{X_k\in I_{|(\alpha
t_{j-1},\alpha t_j]}}\theta_i N_{C_k}(0,\alpha
t_i]\right]=\nonumber\\
&\exp\left(\nu\int_{0}^{\alpha(t_j-
t_{j-1})}\mathrm{E}\left[\exp\left(\sum_{i=j}^{n}\theta_i
N_{C_0}(-\alpha t_{j-1}-s,\alpha
(t_i-t_{j-1})-s]\right)-1\right]\,\mathrm{d}s\right)\label{eq:arsenal}.
\end{align}

\noindent We now show

\begin{equation}\label{Ilimite}
\lim_{\alpha\rightarrow\infty}\frac{1}{\alpha}\log\mathrm{E}\left[\exp\sum_{i=1}^{n}\theta_i
\sum_{X_k\in I_{|(0,\alpha t_i]}}N_{C_k}(0,\alpha t_i]
\right]=\sum_{j=1}^{n}(t_j-t_{j-1})\Lambda\left(\sum_{i=j}^{n}\theta_i\right)
\end{equation}

\noindent for each $(\theta_1,\ldots,\theta_n)\in\mathbb{R}^n$.~We
first notice that by $(\ref{prodotto})$ and $(\ref{eq:arsenal})$
we have

\begin{equation*}
\frac{1}{\alpha}\log\mathrm{E}\left[\exp\sum_{i=1}^{n}\theta_i
\sum_{X_k\in I_{|(0,\alpha t_i]}}N_{C_k}(0,\alpha t_i]
\right]=\sum_{j=1}^{n}(t_j-t_{j-1})J_j(\alpha),
\end{equation*}

\noindent where

\begin{equation}\label{mgf-k2}
J_j(\alpha)= \frac{\nu}{\alpha(t_j-t_{j-1})}
\int_0^{\alpha(t_j-t_{j-1})}\left(\mathrm{E}\left[\exp\sum_{i=j}^{n}\theta_i
N_{C_0}(-\alpha t_{j-1}-s,\alpha
(t_i-t_{j-1})-s]\right]-1\right)\,\mathrm{d}s.
\end{equation}

\noindent Now suppose that
$(\theta_1,\ldots,\theta_n)\in\mathbb{R}^n$ is such that
$\sum_{i=j}^n\theta_i\in\mathcal{D}_{S}$ for each
$j\in\{1,\ldots,n\}$.~Then by Lemma \ref{le:agcmgt} it follows
that there exists $\theta^{*}\in\mathcal{D}_{S}$ such that
$\theta^{*}\geq 0$, $\sum_{i=j}^{n}\theta_i\le\theta^{*}$ for all
$j\in\{1,\ldots,n\}$, and

\[
\sum_{i=j}^{n}\theta_i N_{C_0}(-\alpha t_{j-1}-s,\alpha
(t_i-t_{j-1})-s]\leq\theta^{*}N_{C_0}(\R),\quad\text{a.s..}
\]

\noindent By (\ref{mgf-k2}) and the dominated convergence theorem,
we have

\[
\lim_{\alpha\rightarrow\infty}J_j(\alpha)=\nu\left(\mathrm{E}\left[
\mathrm{e}^{\sum_{i=j}^{n}\theta_{i}S}\right]-1\right).
\]

\noindent Hence we proved $(\ref{Ilimite})$ whenever
$(\theta_1,\ldots,\theta_n)\in\mathbb{R}^{n}$ satisfies
$\sum_{i=j}^{n}\theta_i\in\mathcal{D}_{S}$ for every
$j\in\{1,\ldots,n\}$.~Now suppose that
$(\theta_1,\ldots,\theta_n)\in\mathbb{R}^{n}$ is such that
$\sum_{i=j}^{n}\theta_i\notin\mathcal{D}_{S}$ for some
$j\in\{1,\ldots,n\}$.~We have that $J_j(\alpha)$ is bigger than or
equal to

\begin{align}
&\frac{\nu}{\alpha (t_j-t_{j-1})}\int_0^{\alpha
(t_j-t_{j-1})}\mathrm{E}\left[\exp\left(\sum
_{i=j}^{n}\bold{1}\{\theta _i<0\}\theta_i S+\sum
_{i=j}^{n}\bold{1}\{\theta_i>0\}\theta_i
N_{C_0}[0,\alpha(t_j-t_{j-1})-s]\right)-1\right]\,\mathrm{d}s\nonumber\\
&=\nu\int_0^{1}\mathrm{E}\left[\exp\left(\sum
_{i=j}^{n}\bold{1}\{\theta _i<0\}\theta_i S+\sum
_{i=j}^{n}\bold{1}\{\theta_i>0\}\theta_i N_{C_0}[0,\alpha
(t_j-t_{j-1})(1-z)]\right)-1\right]\,\mathrm{d}z.\nonumber
\end{align}

\noindent The expectation in the latter formula goes to
$\mathrm{E}[\exp(\sum_{i=j}^{n}\theta _i S)-1]$ as
$\alpha\rightarrow\infty$ by the monotone convergence
theorem.~Therefore, by Fatou's lemma we have

\[
\lim_{\alpha\rightarrow\infty}J_j(\alpha)\geq\nu\mathrm{E}\left[\exp
\left(\sum_{i=j}^{n}\theta_i S\right)-1\right]=\infty.
\]

\noindent Thus, since the quantities $J_1(\alpha),\ldots
,J_n(\alpha)$ are bounded below by $-\nu$, we get (\ref{Ilimite})
also in this case.~We now show

\begin{equation}
\label{eq:limIT}
\lim_{\alpha\rightarrow\infty}\frac{1}{\alpha}\log\mathrm{E}\left[\exp\sum_{i=1}^{n}\theta_i
\sum_{X_k\in I_{|[-T,0]}}N_{C_k}(0,\alpha t_i]\right]= 0
\end{equation}

\noindent for all $(\theta_1,\ldots,\theta_n)\in\mathbb{R}^n$ such
that $\sum_{i=j}^{n}\theta_i\in\mathcal{D}_{S}$ for each
$j\in\{1,\ldots,n\}$.~By Lemma \ref{le:agcmgt} we have that there
exists $\theta^{*}\in\mathcal{D}_{S}$ such that $\theta^{*}\geq
0$, $\sum_{i=j}^{n}\theta_i\le\theta^{*}$ for all $j\in\{
1,\ldots,n\}$ and

\[
\theta_{-}\sum_{X_k\in I_{|[-T,0]}}N_{C_k}(\R)
\leq\sum_{i=1}^{n}\theta_i\sum_{X_k\in I_{|[-T,0]}
}N_{C_k}(0,\alpha t_i]\leq\theta^{*}\sum_{X_k\in
I_{|[-T,0]}}N_{C_k}(\R),\quad\text{a.s.},
\]

\noindent where $\theta_{-}\equiv\sum_{i:\theta_i<0}\theta_i$ and
$\theta_{-}\equiv 0$ if $\{i:\theta_i<0\}=\emptyset$.~Therefore,
using again Lemma 6.4 VI in Daley and Vere-Jones (2003), we have

\[
\exp\left(\nu T(\mathrm{E}[\mathrm{e}^{\theta_{-} S}]-1)\right)
\leq \mathrm{E}\left[\exp\sum_{i=1}^{n}\theta_i \sum_{X_k\in
I_{|[-T,0]} }N_{C_k}(0,\alpha t_i]\right]\leq \exp\left(\nu
T(\mathrm{E}[\mathrm{e}^{\theta^{*} S}]-1)\right).
\]

\noindent Equation (\ref{eq:limIT}) follows taking the logarithms
in the above inequalities and passing to the limit.~The conclusion
follows putting together (\ref{eq:nlaplace}), (\ref{Ilimite}) and
(\ref{eq:limIT}).

\vspace{0.2cm}

\noindent {\em Proof of $(\mathrm b)$} {\em and} {\em Proof of
$(\mathrm c)$.} Part {\em $(\mathrm b)$} can be shown using
assumption $(\ref{eq:hypPC})$ and following the lines of the proof
of part {\em $(\mathrm b)$} of Proposition 2.2 in Ganesh, Macci
and Torrisi (2005).~The proof of part {\em $(\mathrm c)$} is
identical to the proof of part {\em $(\mathrm c)$} of Proposition
2.2 in Ganesh, Macci and Torrisi (2005).

\vspace{0.2cm}

\noindent {\em End} {\em of} {\em the} {\em proof} {\em under}
{\em condition} {\em $(\ref{eq:addcond})$.} By claim $(\bold{C})$
and the Dawson-G\"artner theorem,
$\{N_{\bold{X}_{\alpha\cdot,T}}(0,\alpha\cdot]/\alpha\}$ satisfies
the LDP on $D[0,1]$, equipped with the topology of point-wise
convergence, with speed $\alpha$ and good rate function

\[
\tilde{J}(f)=\sup\Bigl\{\sum _{k=1}^n (t_k-t_{k-1})\Lambda^{*}
\Bigl( \frac{f(t_k)-f(t_{k-1})}{t_k-t_{k-1}}\Bigr):\ n\geq 1,
0\leq t_1<\ldots <t_n\leq 1 \Bigr\}.
\]

\noindent The conclusion follows noticing that $\tilde{J}(\cdot)$
coincides with $J(\cdot)$ in (\ref{rfIexpression}), as can be
checked following the same lines as in
the proof of Lemma 5.1.6 in Dembo and Zeitouni (1998).

\vspace{0.2cm}

\noindent {\em Removing} {\em the} {\em additional} {\em condition
$(\ref{eq:addcond})$.} The general case is solved as
follows.~Since $C_k$ is almost surely finite, there exists a
left-most extremal point $Y_k\in C_k$ such that $N_{C_k}(-\infty,
Y_k)=0$ a.s..~Note that, given the immigrants, ${Y_{k}-X_{k}}$ is
an iid sequence.~Therefore, by a classical result on Poisson
processes we have that $\{Y_k\}$ is a Poisson process with
intensity $\nu$.~Viewing ${\bf X}_{t,T}$ as a Poisson cluster
process with cluster centers $Y_k$ and clusters $C_k$, the
conclusion follows by the first part of the proof.
\\
\noindent $\bold{\square}$

\noindent $\bold{Proof\,of\,Theorem\,\ref{SPLDP}}.$ The proof uses
similar steps as in the proof of Theorem
\ref{LDPfortrucated2}.~Here we sketch the main difference.~Assume
the additional condition $N_{C_0}((-\infty,0))=0$ a.s. (the
general case can be treated as in the proof of Theorem
\ref{LDPfortrucated2}).~Define the following subsets of $\R^n$:

\[
A_1\equiv\left\{(\theta_1,\ldots,\theta_n)\in\R^n:\sum_{i=j}^{n}\theta_i\in\D_S\text{
for all $j\in\{1,\ldots,n\}$}\right\}
\]

\noindent and

\[
A_2\equiv\left\{(\theta_1,\ldots,\theta_n)\in\R^n:\sum_{i=j}^{n}\theta_i\notin\D_S\text{
for some $j\in\{1,\ldots,n\}$}\right\}
\]

\noindent We start showing that
for all $n\geq 1$ and $0\leq t_1<\ldots<t_n\leq1$

\begin{equation}\label{eq:PCcor}
\Lambda_{t_1,\ldots ,t_n}(\theta_1,\ldots,\theta_n)=\left \{
\begin{array}{ll}
\sum_{j=1}^{n}(t_j-t_{j-1})\Lambda\left(\sum
_{i=j}^n\theta_i\right) &\ {\rm for}\
(\theta_1,\ldots,\theta_n)\in
A_1\\
+\infty &\ {\rm for}\ (\theta_1,\ldots,\theta_n)\in A_2,
\end{array}
\right.
\end{equation}

\noindent where

\[
\Lambda_{t_1,\ldots ,t_n}(\theta_1,\ldots,\theta_n)\equiv
\lim_{\alpha\rightarrow\infty}\frac{1}{\alpha}\log\mathrm{E}
\left[\exp\left(\sum _{i=1}^n\theta_i N_{\bold{X}}(0,\alpha
t_i]\right)\right]
\]

\noindent and $\Lambda(\cdot)$ is defined in the statement of
Theorem \ref{scalarLDPPC}.~Using the definition of ${\bf X}$ and
the assumption $N_{C_0}((-\infty,0))=0$ a.s., we have

\begin{align}
\mathrm{E}\left[\exp\left(\sum _{i=1}^n\theta_i
N_{\bold{X}}(0,\alpha
t_i]\right)\right]&=\mathrm{E}\left[\exp\sum_{i=1}^{n}\theta_i
\sum_{X_k\in I_{|(0,\alpha t_i]}}N_{C_k}(0,\alpha t_i]
\right]\mathrm{E}\left[\exp\sum_{i=1}^{n}\theta_i \sum_{X_k\in
I_{|[-T,0]}} N_{C_k}(0,\alpha t_i] \right]\times\nonumber\\
&\,\,\,\,\,\,\,\,\,\,\,\,\,\,\,\,\,\,\times
\mathrm{E}\left[\exp\sum_{i=1}^{n}\theta_i \sum_{X_k\in
I_{|(-\infty,-T)}} N_{C_k}(0,\alpha t_i]\right].\nonumber
\end{align}

\noindent As noticed in the proof of Theorem \ref{LDPfortrucated2}
we have

\[
\mathrm{E}\left[\exp\sum_{i=1}^{n}\theta_i N_{\bold{X}_{\alpha
t_i,T}}(0,\alpha t_i]\right]
=\mathrm{E}\left[\exp\sum_{i=1}^{n}\theta_i \sum_{X_k\in
I_{|(0,\alpha t_i]}}N_{C_k}(0,\alpha t_i]
\right]\mathrm{E}\left[\exp\sum_{i=1}^{n}\theta_i \sum_{X_k\in
I_{|[-T,0]}} N_{C_k}(0,\alpha t_i]\right].
\]

\noindent Therefore, by the computations in the proof of Theorem
\ref{LDPfortrucated2}, to prove $(\ref{eq:PCcor})$ we only need to
check that

\begin{equation}\label{eq:finallimitSP}
\lim_{\alpha\to \infty}\frac{1}{\alpha}\log
\mathrm{E}\left[\exp\sum_{i=1}^{n}\theta_i \sum_{X_k\in
I_{|(-\infty,-T)}} N_{C_k}(0,\alpha t_i]\right]=0,\quad\text{for
all $(\theta_1,\ldots,\theta_n)\in A_1$.}
\end{equation}

\noindent Since $\{(X_i,C_i):X_i\in I_{|(-\infty,-T)}\}$ is an
independently marked Poisson process, by Lemma 6.4.VI in Daley and
Vere-Jones (2003) we have

\begin{align}
\mathrm{E}\left[\exp\sum_{i=1}^{n}\theta_i \sum_{X_k\in
I_{|(-\infty,-T)}} N_{C_k}(0,\alpha t_i]\right]&= \mathrm{E}\left[
\exp\sum_{X_k\in I_{|(-\infty,-T)}}\sum_{i=1}^{n}\theta_i
N_{C_k-X_k}(-X_k,\alpha
t_i-X_k]\right]\nonumber\\
&=\exp\left(\nu\int_{T}^{\infty}\mathrm{E}\left[\mathrm{e}^{\sum_{i=1}^{n}\theta_i
N_{C_0}(x,\alpha t_i+x]}-1\right]\,\mathrm{d}x\right)\nonumber
\end{align}

\noindent Take $(\theta_1,\ldots,\theta_n)\in A_1$.~By Lemma
\ref{le:agcmgt} we have that there exists
$\theta^{*}\in\mathcal{D}_{S}$ such that $\theta^{*}\geq 0$,
$\sum_{i=j}^{n}\theta_i\le\theta^{*}$ for all $j\in\{
1,\ldots,n\}$ and

\[
\theta_{-}N_{C_0}(\R)\leq\sum_{i=1}^{n}\theta_i N_{C_0}(x,\alpha
t_i+x]\leq\theta^{*}N_{C_0}(\R),\quad\text{a.s.}
\]

\noindent where $\theta_{-}\equiv\sum_{i:\theta_i<0}\theta_i$ and
$\theta_{-}\equiv 0$ if $\{i:\theta_i<0\}=\emptyset$.~Thus,

\[
\mathrm{e}^{\sum_{i=1}^{n}\theta_i N_{C_0}(x,\alpha
t_i+x]}-1\leq(\mathrm{e}^{\theta^{*}N_{C_0}(\R)}-1)\bold{1}\{x\leq
L\},\quad\text{for all $x\geq T$}
\]

\noindent and

\[
\mathrm{e}^{\sum_{i=1}^{n}\theta_i N_{C_0}(x,\alpha
t_i+x]}-1\geq(\mathrm{e}^{\theta_{-}N_{C_0}(\R)}-1)\bold{1}\{x\leq
L\},\quad\text{for all $x\geq T$}
\]

\noindent The limit $(\ref{eq:finallimitSP})$ follows by
assumption $(\ref{eq:radCo})$ noticing that the above relations
yield, for all $(\theta_1,\ldots,\theta_n)\in A_1$:

\[
\mathrm{E}\left[\exp\sum_{i=1}^{n}\theta_i \sum_{X_k\in
I_{|(-\infty,-T)}} N_{C_k}(0,\alpha t_i]\right]
\leq\exp(\nu\mathrm{E}[L(\mathrm{e}^{\theta^* S}-1)])
\]

\noindent and

\[
\mathrm{E}\left[\exp\sum_{i=1}^{n}\theta_i \sum_{X_k\in
I_{|(-\infty,-T)}} N_{C_k}(0,\alpha
t_i]\right]\geq\exp(\nu\mathrm{E}[(\mathrm{e}^{\theta_-
S}-1)(L-T)\bold{1}\{L\geq T\}]).
\]

\noindent Now since $\D_S$ is open, the claim follows by applying
first the G{\"a}rtner-Ellis theorem in $\R^n$ to get the LDP for
the finite-dimensional distributions, and then the
Dawson-G{\"a}rtner theorem to have the LDP for the process (argue
as in the proof of Theorem \ref{LDPfortrucated2} for the remaining
steps).\\
\noindent$\bold{\square}$

\subsection{Sample path large deviations in the topology of uniform convergence}
In the applications, one usually derives LDPs for continuous
functions of sample paths of stochastic processes by using the
contraction principle.~Since the topology of uniform convergence
is finer than the topology of point-wise convergence, it has a
larger class of continuous functions.~Thus, it is of interest to
understand if $\{N_{\bold{X}}(0,\alpha\cdot]/\alpha\}$ satisfies a
LDP on $D[0,1]$ equipped with the topology of uniform
convergence.~In this section we give an answer to this question
assuming that the tails of $S$ decay super-exponentially.

\begin{Theorem}\label{SPLDPunif}
Assume

\begin{equation}\label{superexp}
\mathrm{E}[\mathrm{e}^{\theta S}]<\infty\quad\text{for each
$\theta\in\R$}
\end{equation}

\noindent and

\begin{equation}\label{eq:radCosup_ex}
\mathrm{E}[L\mathrm{e}^{\theta S}]<\infty\quad\text{for each
$\theta\in\R$.}
\end{equation}

\noindent Then $\{\frac{N_{\bold{X}}(0,\alpha \cdot]}{\alpha}\}$
satisfies a LDP on $D[0,1]$, equipped with the topology of uniform
convergence, with speed $\alpha$ and good rate function
$(\ref{rfIexpression})$.
\end{Theorem}

\noindent In this section, without loss of generality we assume
that the points of $I$ are $\{X_i\}_{i\in\mathbb{Z}^*}$, where
$\mathbb{Z}^*=\mathbb{Z}\backslash\{0\}$, $X_i<X_{i+1}$, and we
set $X_0=0$.~As usual, we denote by ${\bf X}_{t,T}$ the
non-stationary Poisson cluster process with immigrant process
empty on $(-\infty,-T)\cup (t+T,\infty)$, where $T>0$ is a fixed
constant, and by $C_0$ the cluster generated by an immigrant at
$0$.

Before proving Theorem \ref{SPLDPunif} we show that the same LDP
holds for ${\bf X}_{t,T}$, under a weaker condition.

\begin{Theorem}\label{LDPfortrucated2unif}
Assume $(\ref{superexp})$.~Then
$\{N_{\bold{X}_{\alpha\cdot,T}}(0,\alpha\cdot]/\alpha\}$ satisfies
a LDP on $D[0,1]$, equipped with the topology of uniform
convergence, with speed $\alpha$ and good rate function
$(\ref{rfIexpression})$.
\end{Theorem}

\noindent To prove Theorem \ref{LDPfortrucated2unif} above we use
the following Lemma \ref{le:agcmgt1}, whose proof is omitted since
it is similar to the proof of Lemma 3.3 in Ganesh, Macci and
Torrisi (2005).~Let $\{S_k\}_{k\in\mathbb{Z}}$ be the iid sequence
of random variables (distributed as $S$) defined by $S_k = N_{C_k}
(\R)$.

\begin{Lemma}\label{le:agcmgt1}
Assume $(\ref{superexp})$, $N_{C_0}(-\infty,0)=0$ a.s., and define

\begin{equation*}
A_{n}=\sum_{k=0}^{n-1}(S_k-N_{C_k-X_k}(0,X_k]),\quad\text{$n\geq
1$.}
\end{equation*}

\noindent It holds

\begin{equation*}
\lim_{n\to\infty}\frac{1}{n}\log P(A_n\geq
n\delta)=-\infty\quad\text{for each $\delta>0$.}
\end{equation*}

\end{Lemma}

\noindent $\bold{Proof\,of\,Theorem\,\ref{LDPfortrucated2unif}}.$
We prove the theorem assuming that $N_{C_0}(-\infty,0)=0$
a.s..~The general case is solved as in the proof of Theorem
\ref{LDPfortrucated2}.~As usual denote by $I_{|A}$ the restriction
of $I$ on the Borel set $A\subset\mathbb{R}$.~Define

\begin{equation*}
C(t)=\sum_{X_k\in I_{|(0,t]}}S_k,\quad\text{$t>0$.}
\end{equation*}

\noindent We prove that
$\{N_{\bold{X}_{\alpha\cdot,T}}(0,\alpha\cdot]/\alpha\}$ and
$\{C(\alpha\cdot)/\alpha\}$ are exponentially equivalent (see, for
instance, Definition 4.2.10 in the book of Dembo and Zeitouni,
(1998)) with respect to the topology of uniform
convergence.~Therefore the conclusion follows by a well-known
result on sample path large deviations, with respect to the
uniform topology, of compound Poisson processes (see, for
instance, Borovkov (1967); see also de Acosta (1994) and the
references cited therein) and Theorem 4.2.13 in Dembo and Zeitouni
(1998).~Define

\begin{equation*}
C_{T}(t)=\sum_{X_k\in I_{|[-T,t]}}S_k,\quad\text{$t>0$.}
\end{equation*}

\noindent Using Chernoff bound and condition $(\ref{superexp})$
can be easily realized that the processes
$\{C(\alpha\cdot)/\alpha\}$ and $\{C_{T}(\alpha\cdot)/\alpha\}$
are exponentially equivalent with respect to the topology of
uniform convergence.~Therefore, it suffices to show that
$\{C_{T}(\alpha\cdot)/\alpha\}$ and
$\{N_{\bold{X}_{\alpha\cdot,T}}(0,\alpha\cdot]/\alpha\}$ are
exponentially equivalent with respect to the topology of uniform
convergence.~Note that the assumption $N_{C_0}(-\infty,0)=0$ a.s.
gives

\begin{equation*}
N_{\bold{X}_{t,T}}(0,t]=\sum_{X_k\in
I_{|(0,t]}}N_{C_k}(0,t]+\sum_{X_k\in
I_{|[-T,0]}}N_{C_k}(0,t]\quad\text{$t>0$, a.s..}
\end{equation*}

\noindent Therefore, we need to show that

\begin{equation}\label{EEcondition}
\lim_{\alpha \to \infty} \frac{1}{\alpha}\log P(M_{\alpha}>\delta)
=-\infty,\quad\text{for any $\delta>0$,}
\end{equation}

\noindent where

\[
M_{\alpha}=\frac{1}{\alpha}\sup_{t\in [0,1]}\Big|C_T(\alpha
t)-\sum_{X_k\in I_{|(0,\alpha t]}}N_{C_k}(0,\alpha t]-\sum_{X_k\in
I_{|[-T,0]}}N_{C_k}(0,\alpha t]\Big|.
\]

\noindent Since

\[
M_\alpha\leq M_{\alpha}^{(1)}+M_{\alpha}^{(2)},\quad\text{a.s.,}
\]

\noindent where

\begin{equation*}
M_{\alpha}^{(1)}=\frac{1}{\alpha}\sum_{X_k\in
I_{|[-T,0]}}S_k\quad\text{and}\quad
M_{\alpha}^{(2)}=\frac{1}{\alpha}\sup_{t\in [0,1]}\sum_{X_k\in
I_{|(0,\alpha t]}}(S_k-N_{C_k}(0,\alpha t]),
\end{equation*}

\noindent the limit $(\ref{EEcondition})$ follows if we prove

\begin{equation}\label{EEcondition1}
\lim_{\alpha \to \infty} \frac{1}{\alpha}\log
P(M_{\alpha}^{(1)}>\delta/2)=-\infty,\quad\text{for any
$\delta>0$}
\end{equation}

\noindent and

\begin{equation}\label{EEcondition2}
\lim_{\alpha \to \infty} \frac{1}{\alpha}\log
P(M_{\alpha}^{(2)}>\delta/2)=-\infty,\quad\text{for any
$\delta>0$.}
\end{equation}

\noindent The limit $(\ref{EEcondition1})$ easily follows by the
Chernoff bound and condition $(\ref{superexp})$.~It remains to
show $(\ref{EEcondition2})$.~Since the random function $t\mapsto
N_{C_k}(0,\alpha t]$ is non-decreasing, it is clear that the
supremum over $t$ is attained at one of the points $X_n$, $n\geq
1$.~Thus

\begin{equation*}
M_{\alpha}^{(2)}=\frac{1}{\alpha}\max_{n\geq 1:X_n\le
\alpha}\sum_{k=1}^{n}(S_k-N_{C_k}(0,X_n]).
\end{equation*}

\noindent Note that

\begin{equation*}
M_{\alpha}^{(2)}\leq\tilde{M}_{\alpha}\quad\text{where
$\tilde{M}_{\alpha}=\frac{1}{\alpha}\max_{n\geq 1:X_n\le
\alpha}\sum_{k=1}^{n}(S_k-N_{C_k-X_k}(0,X_n-X_k])$ a.s..}
\end{equation*}

\noindent Therefore $(\ref{EEcondition2})$ follows if we show

\begin{equation}\label{EEcondition3}
\lim_{\alpha\to\infty}\frac{1}{\alpha}\log
P(\tilde{M}_{\alpha}>\delta/2)=-\infty,\quad\text{for any
$\delta>0$.}
\end{equation}

\noindent Since $X_n$, $n\geq 1$, is the sum of $n$ exponential
random variables with mean $1/\nu$, using Chernoff bound and
taking the logarithm, we have that, for all $\eta>0$ and all
integers $K>\nu$,

\begin{equation} \label{bound1}
\frac{1}{\alpha}\log P(X_{K[\alpha]}<\alpha)\leq\eta+
\frac{K[\alpha]}{\alpha}\log\Bigl(\frac{\nu}{\nu+\eta}\Bigr).
\end{equation}

\noindent Here the symbol $[\alpha]$ denotes the integer part of
$\alpha$.~Next, observe that using the union bound we get

\begin{equation*}
P(\tilde{M}_{\alpha}>\delta/2,X_{K[\alpha]}\ge\alpha)\le K
[\alpha]\max_{1\le n\le K[\alpha]}P\Bigl(\sum_{k=1}^n
(S_k-N_{C_k-X_k}(0,X_{n}-X_{k}])\ge\alpha\delta/2\Bigr),
\end{equation*}

\noindent Now we remark that, for $n\geq 1$, $(X_n-X_1,\ldots
,X_n-X_{n-1})$ and $(X_{n-1},\ldots,X_1)$ have the same joint
distribution.~Moreover, given $I$, the centered processes
$C_k-X_k$ are iid and independent of the $\{X_k\}$.~Hence, letting
$A_n$ denote the random variable defined in the statement of Lemma
\ref{le:agcmgt1}, we have

\begin{equation*}
P(\tilde{M}_{\alpha}>\delta/2,X_{K[\alpha]}\ge\alpha)\leq
K[\alpha]\max_{1\le n\le K[\alpha]}P(A_{n}\geq\alpha\delta/2).
\end{equation*}

\noindent The random variables $A_n$ are increasing in $n$,
therefore

\[
P(\tilde{M}_{\alpha}>\delta,X_{K[\alpha]}\ge\alpha)\le K[\alpha]P(
A_{K[\alpha]}\ge\alpha\delta/2),
\]

\noindent and by Lemma \ref{le:agcmgt1} we have

\begin{equation}\label{bound2}
\lim_{\alpha \to \infty} \frac{1}{\alpha}\log
P(\tilde{M}_{\alpha}>\delta/2,X_{K[\alpha]}\ge\alpha)=-\infty.
\end{equation}

\noindent Now note that

\begin{equation*}
P(\tilde{M}_{\alpha}>\delta/2)\leq
P(\tilde{M}_{\alpha}>\delta/2,X_{K[\alpha]}\geq\alpha)+
P(X_{K[\alpha]}<\alpha),
\end{equation*}

\noindent for arbitrary $K>\nu$.~Hence by (\ref{bound1})
and(\ref{bound2}) we have

$$
\limsup_{\alpha\to\infty}\frac{1}{\alpha}\log
P(\tilde{M}_{\alpha}>\delta/2)\leq\inf_{\eta >0}\left(\eta+K\log
\Bigl(\frac{\nu}{\nu+\eta}\Bigr)\right)= K-\nu-K\log\frac{K}{\nu}.
$$

\noindent Then we obtain (\ref{EEcondition3}) by letting $K$ tend
to $\infty$.
\\
\noindent$\square$

\noindent $\bold{Proof\,of\,Theorem\,\ref{SPLDPunif}}.$ Throughout
the proof we assume $N_{C_0}((-\infty,0))=0$ a.s..~The general
case is solved as in the proof of Theorem
\ref{LDPfortrucated2}.~Let $\{C_T(t)\}$ be the process defined in
the proof of Theorem \ref{LDPfortrucated2unif}.~The claim follows
if we show that $\{C_{T}(\alpha\cdot)/\alpha\}$ and
$\{N_{\bold{X}}(0,\alpha\cdot]/\alpha\}$ are exponentially
equivalent with respect to the topology of uniform
convergence.~Note that the assumption $N_{C_0}(-\infty,0)=0$ a.s.
implies

\begin{equation*}
N_{\bold{X}}(0,t]=N_{\bold{X}_{t,T}}(0,t]+\sum_{X_k\in
I_{|(-\infty,-T)}}N_{C_k}(0,t],\quad\text{$t>0$ a.s..}
\end{equation*}

\noindent Therefore, since we already proved that
$\{C_{T}(\alpha\cdot)/\alpha\}$ and
$\{N_{\bold{X}_{\alpha\cdot,T}}(0,\alpha\cdot]/\alpha\}$ are
exponentially equivalent with respect to the uniform topology (see
the proof of Theorem \ref{LDPfortrucated2unif}), the claim follows
if we prove that

\begin{equation}\label{EEEcondition}
\lim_{\alpha \to \infty} \frac{1}{\alpha}\log P\left(\sum_{X_k\in
I_{|(-\infty,-T)}}N_{C_k}(0,\alpha]>\alpha\delta\right)
=-\infty,\quad\text{for any $\delta>0$.}
\end{equation}

\noindent Using the Chernoff bound we have, for all $\theta>0$,

\begin{align}
P\left(\sum_{X_k\in
I_{|(-\infty,-T)}}N_{C_k}(0,\alpha]>\alpha\delta\right)&\leq\mathrm{e}^{-\alpha\theta\delta}
\mathrm{E}\left[\exp\sum_{X_k\in I_{|(-\infty,-T)}}\theta
N_{C_k}(0,\alpha]\right]\nonumber\\
&=\mathrm{e}^{-\alpha\theta\delta}
\exp\left(\nu\int_{T}^{\infty}\mathrm{E}[\mathrm{e}^{\theta
N_{C_0}(x,\alpha+x]}-1]\,\mathrm{d}x\right)\nonumber\\
&\leq\mathrm{e}^{-\alpha\theta\delta}\exp\left(\nu\mathrm{E}[(\mathrm{e}^{\theta
S}-1)L]\right).\nonumber
\end{align}

\noindent Taking the logarithm, dividing by $\alpha$, letting
$\alpha$ tend to $\infty$ and using assumption
$(\ref{eq:radCosup_ex})$ we get

\[
\limsup_{\alpha\to\infty}\frac{1}{\alpha}\log P\left(\sum_{X_k\in
I_{|(-\infty,-T)}}N_{C_k}(0,\alpha]>\alpha\delta\right)\leq-\theta\delta,\quad\text{for
all $\theta>0$.}
\]

\noindent Relation $(\ref{EEEcondition})$ follows letting $\theta$
tend to infinity in the above inequality.
\\
\noindent$\bold{\square}$

\section{Large deviations of spatial Poisson cluster processes}
\label{sec:LDPSPAT}

\subsection{The large deviations principle}
A spatial Poisson cluster process ${\bf X}$ is a Poisson cluster
process in $\R^d$, where $ d\geq 1$ is an integer.~The clusters
centers are the points $\{X_i\}$ of a homogeneous Poisson process
$I\subset\R^d$ with intensity $\nu\in (0,\infty)$.~Each immigrant
$X_i\in I$ generates a cluster $C_i=C_{X_i}$, which is a finite
point process.~Given $I$, the centered clusters $\{C_{X_i}-X_i\}$
are iid and independent of $I$.~${\bf X}$ is the union of all
clusters.~As in dimension $1$, we denote by $S$ the number of
points in a cluster, with a little abuse of notation by $C_0$ the
cluster generated by a point at $0$, and by $L$ the radius of
$C_0$.~Moreover, we denote by $N_{\bold{X}}(b(0,r))$ the number of
points of $\bold{X}$ in the ball $b(0,r)$, and by

\[
\omega_{d}(r)=\frac{r^d\pi^{d/2}}{\Gamma(1+d/2)}
\]

\noindent the volume of $b(0,r)$.~The following LDP holds:

\begin{Theorem}\label{scalarLDPspatPC}
Assume $(\ref{eq:hypPC})$ and

\begin{equation}\label{eq:radCoSp}
\mathrm{E}[L^d\mathrm{e}^{\theta S}]<\infty,\quad\text{for all
$\theta\in\mathcal{D}_S^{\circ}$.}
\end{equation}

\noindent Then $\{N_{\bold{X}}(b(0,r))/\omega_d(r)\}$ satisfies a
LDP on $\mathbb{R}$ with speed $\omega_d(r)$ and good rate
function $(\ref{rfPC})$.
\end{Theorem}

\noindent Before proving Theorem \ref{scalarLDPspatPC}, we show
that the same LDP holds for the non-stationary Poisson cluster
process $\bold{X}_{r,R}$ with immigrant process empty in
$\mathbb{R}^{d}\setminus b(0,R+r)$.~As usual this LDP holds under
a weaker condition.

\begin{Theorem}\label{scalarLDPspatPCtru}
Assume $(\ref{eq:hypPC})$.~Then
$\{N_{\bold{X}_{r,R}}(b(0,r))/\omega_d(r)\}$ satisfies a LDP on
$\mathbb{R}$ with speed $\omega_d(r)$ and good rate function
$(\ref{rfPC})$.
\end{Theorem}

\noindent $\bold{Proof\,of\,Theorem\, \ref{scalarLDPspatPCtru}}$
The proof is similar to that one for the non-stationary Poisson
cluster process on the line.~Here we just sketch the main
differences.~As in the proof of Theorem \ref{LDPfortrucatedPC},
the claim follows by the G{\"a}rtner-Ellis theorem.~Indeed,
letting $I_{| b(0,r)}$ denote the point process of immigrants in
$b(0,r)$, and $I_{|b(0,R+r)\setminus b(0,r) } $ the point process
of immigrants in $b(0,R+r)\setminus b(0,r)$ we have, for each
$\theta\in\mathbb{R}$,

\begin{equation*}
\mathrm{E}\left[\mathrm{e}^{\theta
N_{\bold{X}_{r,R}}(b(0,r))}\right]=
\mathrm{E}\left[\mathrm{e}^{\theta\sum_{X_i\in I_{|
b(0,r)}}N_{C_i}(b(0,r))}\right]\mathrm{E}\left[\mathrm{e}^{\theta\sum_{X_i\in
I_{|b(0,R+r)\setminus b(0,r) } } N_{C_i}(b(0,r))}\right].\nonumber
\end{equation*}

\noindent As usual, with a little abuse of notation denote by
$C_0$ the cluster generated by an immigrant at $0$.~It holds:

\[
\mathrm{E}\left[\mathrm{e}^{\theta\sum_{X_i\in I_{|
b(0,r)}}N_{C_i}(b(0,r))}\right]=\exp\left(\nu\int_{b(0,r)}\mathrm{E}\left[\mathrm{e}^{\theta
N_{C_0}(b(-x,r))}-1\right]\,\mathrm{d}x\right),
\]

\noindent for each $\theta\in\mathbb{R}$;

\[
\mathrm{E}\left[\mathrm{e}^{\theta\sum_{X_i\in
I_{|b(0,R+r)\setminus b(0,r) }
}N_{C_i}(b(0,r))}\right]\leq\exp\left(\nu ( \omega_{d}(R+r) -
\omega_{d}(r))\mathrm{E}\left[\mathrm{e}^{\theta
S}-1\right]\right),
\]

\noindent for $\theta \in [0,\infty) \cap \mathcal D_S$;

\[
1\geq\mathrm{E}\left[\mathrm{e}^{\theta\sum_{X_i\in
I_{|b(0,R+r)\setminus b(0,r)
}}N_{C_i}(b(0,r))}\right]\geq\exp\left(\nu( \omega_{d}(R+r) -
\omega_{d}(r))\mathrm{E}\left[\mathrm{e}^{\theta
S}-1\right]\right),
\]

\noindent for $\theta\leq 0$.~Therefore,

\[
\begin{aligned}
\lim_{r\rightarrow\infty}\frac{1}{\omega_{d}(r)}\log\mathrm{E}\left[\mathrm{e}^{\theta\sum_{X_i\in
I_{|b(0,r)}}N_{C_i}(b(0,r))}\right]&=\lim_{r\rightarrow\infty}\frac{\nu}{\omega_d(r)}\int_{b(0,r)}
\mathrm{E}[\mathrm{e}^{\theta
N_{C_0}(b(-x,r))}-1]\,\mathrm{d}x\nonumber\\
&=\frac{\nu}{\omega_d(1)}\lim_{r\rightarrow\infty}\int_{b(0,1)}\mathrm{E}[\mathrm{e}^{\theta
N_{C_0}(b(-ry,r))}-1]\,\mathrm{d}y\nonumber\\
&=\nu\mathrm{E}[\mathrm{e}^{\theta S}-1],\qquad\text{for each
$\theta\in\mathbb{R}$,}\nonumber
\end{aligned}
\]

\noindent and, since $\lim_{ r\to \infty}
\omega_{d}(R+r)/\omega_{d}(r)=1$, for each $\theta\in\mathcal
D_S$,

\[
\lim_{r\rightarrow\infty}\frac{1}{\omega_d(r)}\log\mathrm{E}\left[\mathrm{e}^{\theta\sum_{X_i\in
I_{|b(0,R+r)\setminus b(0,r)}}N_{C_i}(b(0,r))}\right]=0.
\]

\noindent The rest of the proof is exactly as in the
one-dimensional case.\\
\noindent$\bold{\square}$

\noindent $\bold{Proof\,of\,Theorem\, \ref{scalarLDPspatPC}}$ The
proof is similar to that one of Theorem \ref{scalarLDPspatPCtru}
and is again based on the G{\"a}rtner-Ellis theorem.~We start
showing that

\begin{equation}\label{eq:finallimitspat}
\lim_{r\to
\infty}\frac{1}{\omega_d(r)}\log\mathrm{E}\left[\mathrm{e}^{\theta
(N_{{\bf X}}(b(0,r))-N_{{\bf X}
_{r,R}}(b(0,r)))}\right]=0,\quad\text{for all
$\theta\in\D_S^{\circ}$.}
\end{equation}

\noindent This relation is easily verified for $\theta\leq
0$.~Thus we only check $(\ref{eq:finallimitspat})$ for
$\theta\in\D_S^{\circ}\cap (0,\infty)$.~We have, for all
$\theta\in\D_S^{\circ}\cap (0,\infty)$,

\begin{align}
\mathrm{E}\left[\mathrm{e}^{\theta (N_{{\bf X}}(b(0,r))-N_{{\bf X}
_{r,R}}(b(0,r)))}\right]&=\mathrm{E}\left[\mathrm{e}^{\theta\sum_{X_i\in
I_{|\R^d\setminus b(0,r+R)}}N_{C_i}(b(0,r))}\right]\nonumber\\
&=\exp\left(\nu\int_{\R^d\setminus
b(0,r+R)}\mathrm{E}\left[\mathrm{e}^{\theta
N_{C_0}(b(-x,r))}-1\right]\,\mathrm{d}x\right)\nonumber.
\end{align}

\noindent Now notice that since $\theta>0$ we have

\[
\mathrm{e}^{\theta N_{C_0}(b(-x,r))}-1\leq(\mathrm{e}^{\theta
N_{C_0}(\R^d)}-1)\bold{1}\{\|x\|\leq L+r\},\quad\text{for all
$x\in\R^d$.}
\]

\noindent The limit $(\ref{eq:finallimitspat})$ follows by
assumption $(\ref{eq:radCoSp})$ noticing that the above relations
yield, for all $\theta\in\D_S^{\circ}\cap (0,\infty)$, $r>0$,

\[
\mathrm{E}\left[\mathrm{e}^{\theta\sum_{X_i\in I_{|\R^d\setminus
b(0,r+R)}}N_{C_i}(b(0,r))}\right]\leq\exp((\nu\pi^{d/2}/\Gamma(1+d/2))\mathrm{E}[(L+r)^d(\mathrm{e}^{\theta
S}-1)]).
\]

\noindent Now notice that

\[
\mathrm{E}\left[\mathrm{e}^{\theta N_{{\bf
X}}(b(0,r))}\right]=\mathrm{E}\left[\mathrm{e}^{\theta N_{{\bf X}
_{r,R}}(b(0,r))}\right]\mathrm{E}\left[\mathrm{e}^{\theta (N_{{\bf
X}}(b(0,r))-N_{{\bf X}_{r,R}}(b(0,r)))}\right],\quad\text{for all
$\theta\in\R$, $r>0$.}
\]

\noindent Therefore, if $\D_S=\D_S^{\circ}$ then the claim is a
consequence of the computation of the log-Laplace limit of
$\{N_{{\bf X}_{r,R}}(b(0,r))\}$ in the proof of Theorem
\ref{scalarLDPspatPCtru}, the G{\"a}rtner-Ellis theorem and
assumption (\ref{eq:hypPC}).~It remains to deal with the case
$\D_S\neq\D_S^{\circ}$.~Arguing exactly as in the proof of Theorem
\ref{scalarLDPPC} it can be proved that for any sequence
$\{r_n\}_{n\geq 1}\subset (0,\infty)$ diverging to $+\infty$, as
$n\to\infty$, there exists a subsequence $\{q_n\}\subseteq\{r_n\}$
such that

\[
\limsup_{n\to\infty}\frac{1}{q_n}\log P(N_{\bf
X}(b(0,q_n))/\omega_d(q_{n})\in F)\leq -\inf_{x\in
F}\Lambda^{*}(x),\quad\text{for all closed sets $F$}
\]

\noindent and

\[
\liminf_{n\to\infty}\frac{1}{q_n}\log P(N_{\bf
X}(b(0,q_n))/\omega_d(q_n)\in G)\geq-\inf_{x\in
G}\Lambda^{*}(x),\quad\text{for all open sets $G$,}
\]

\noindent where $\Lambda^*$ is defined by $(\ref{rfPC})$.~Then the
large deviations upper and lower bounds hold for any sequence
$\{r_n\}$ and the claim follows.\\
\noindent$\bold{\square}$

\subsection{The asymptotic behavior of the void probability function and the empty space function}
\label{sec:void} Apart some specific cases, the void probability
function $v(r)=P(N_{\bold{X}}(b(0,r))=0)$, $r>0$, of a spatial
Poisson cluster process is not known in closed form.~Comparing
$\bold{X}$ with the immigrant process $I$ we easily obtain

\begin{equation}\label{voidin}
v(r)\leq P(N_I(b(0,r))=0)=\mathrm{e}^{-\nu\omega_d(r)},\quad r>0.
\end{equation}

\noindent A more precise information on the asymptotic behavior of
$v(\cdot)$, as $r\rightarrow\infty$, is provided by the following
proposition:

\begin{Proposition}\label{cor:void}
Assume $\mathrm{E}[L^d]<\infty$.~Then

\[
\lim_{r\rightarrow\infty}\frac{1}{\omega_d(r)}\log v(r)=-\nu.
\]

\end{Proposition}

\noindent $\bold{Proof}$ Note that

\begin{align}
v(r)&=P(N_{C_i}(b(0,r))=0,\,\text{for all $X_i\in I$})\nonumber\\
&=\mathrm{E}\left[\bold{1}\{N_I(b(0,r))=0\}\prod_{X_i\in
I_{|\R^d\setminus
b(0,r)}}\bold{1}\{N_{C_i}(b(0,r))=0\}\right]\nonumber\\
&=\mathrm{e}^{-\nu\omega_d(r)}\mathrm{E}\left[\prod_{X_i\in
I_{|\R^d\setminus
b(0,r)}}\bold{1}\{N_{C_i}(b(0,r))=0\}\right]\nonumber\\
&=\mathrm{e}^{-\nu\omega_d(r)}\exp\left(-\nu\int_{\R^d\setminus
b(0,r)}P(N_{C_0}(b(-x,r))>0)\,\mathrm{d}x\right)\label{cor:cor}
\end{align}

\noindent where in $(\ref{cor:cor})$ we used Lemma 6.4.VI in Daley
and Vere-Jones (2003).~Thus the claim follows if we prove

\[
\lim_{r\to\infty}\frac{1}{\omega_d(r)}\int_{\R^d\setminus
b(0,r)}P(N_{C_0}(b(-x,r))>0)\,\mathrm{d}x=0.
\]

\noindent For this note that

\[
\bold{1}\{N_{C_0}(b(0,r))>0\}\leq\bold{1}\{\|x\|-r\leq
L\},\quad\text{for all $x\in\R^d$, $r>0$}.
\]

\noindent Therefore,

\[
\frac{1}{\omega_d(r)}\int_{\R^d\setminus
b(0,r)}P(N_{C_0}(b(-x,r))>0)\,\mathrm{d}x\leq\mathrm{E}[(1+L/r)^d-1],
\]

\noindent and the right-hand side in the above inequality goes to
$0$ as $r\to\infty$ by the dominated convergence theorem (note
that $\mathrm{E}[L^d]<\infty$ by assumption).
\\
\noindent$\bold{\square}$

\noindent In spatial statistics, a widely used summary statistic
is the so-called empty space function, which is the distribution
function of the distance from the origin to the nearest point in
$\bold{X}$ (see, for instance, M{\o}ller and Waagepetersen
(2004)), that is

\[
e(r)=1-v(r),\quad\text{$r>0$.}
\]

\noindent Apart some specific cases, the empty space function of
Poisson cluster processes  seems to be intractable.~Next Corollary
\ref{cor:empty} concerns the asymptotic behavior of $e(r)$, as
$r\rightarrow\infty$.

\begin{Corollary}\label{cor:empty}
Under the assumption of Proposition \ref{cor:void} it holds

\[
\lim_{r\rightarrow\infty}\frac{1}{\omega_d(r)}\log\log
e(r)^{-1}=-\nu.
\]

\end{Corollary}

\noindent $\bold{Proof}$ The proof is an easy consequence of
Proposition \ref{cor:void}.~By the upper bound $(\ref{voidin})$ we
obtain

\begin{equation*}
\limsup_{r\rightarrow\infty}\frac{1}{\omega_d(r)}\log\log
e(r)^{-1}\leq\lim_{r\rightarrow\infty}\frac{1}{\omega_d(r)}\log\log
(1-\mathrm{e}^{-\nu\omega_d(r)})^{-1}=-\nu.
\end{equation*}

\noindent To get the matching lower bound we note that the
inequality $\log(1-x)\leq -x$, $x\in [0,1)$, gives $\log
e(r)^{-1}\geq v(r)$, $r>0$, and therefore by Proposition
\ref{cor:void} we get

\[
\liminf_{r\rightarrow\infty}\frac{1}{\omega_d(r)}\log\log
e(r)^{-1}\geq\liminf_{r\rightarrow\infty}\frac{1}{\omega_d(r)}\log
v(r)=-\nu.
\]

\noindent$\bold{\square}$

\subsection{Spatial Hawkes processes}
Spatial Hawkes processes have been introduced in Daley and
Vere-Jones (2003).~Br{\'e}maud, Massouli{\'e} and Ridolfi (2005)
considered spatial Hawkes processes with random fertility rate and
not necessarily Poisson immigrants, and computed the Bartlett
spectrum; the reader is directed to Daley and Vere-Jones (2003)
for the notion of Bartlett spectrum.~M{\o}ller and Torrisi (2005)
derived the pair correlation function of spatial Hawkes processes;
we refer the reader to M{\o}ller and Waagepetersen (2004) for the
notion of pair correlation function.

For the sake of completeness, we briefly recall the definition of
spatial Hawkes process.~A spatial Hawkes process is a Poisson
cluster process $\bold{X}\subset\mathbb{R}^d$ where $d\geq 1$ is
an integer.~The clusters centers are the points $\{X_i\}$ of a
homogeneous Poisson process $I\subset\mathbb{R}^d$ with intensity
$\nu\in (0,\infty)$.~Each immigrant $X_i\in I$ generates a cluster
$C_i=C_{X_i}$ which is formed by the points of generations
$n=0,1,\ldots$ with the following branching structure: the
immigrant $X_i\in I$ is said to be of zero-th generation.~Given
generations $0,1,\ldots,n$ in $C_i$, each point $Y\in C_i$ of
generation $n$ generates a Poisson process on $\mathbb{R}^d$ of
offspring of generation $n+1$ with intensity function
$h(\cdot-Y)$.~Here $h:\mathbb{R}^d\rightarrow [0,\infty)$ is a
non-negative Borel function.~In the model it is assumed that,
given the immigrants, the centered clusters $\{C_{i}-X_{i}\}$ are
iid, and independent of $I$.~By definition the spatial Hawkes
process is $\bold{X}\equiv\bigcup_{i}C_i$.~As in the
one-dimensional case, it is assumed

\begin{equation}\label{stability2}
0<\mu\equiv\int_{\mathbb{R}^d}h(\xi)\,\mathrm{d}\xi<1.
\end{equation}

\noindent This assumption guarantees that the number of points in
a cluster has a finite mean equal to $1/(1-\mu)$, excludes the
trivial case where there are no offspring, and ensures that
$\bold{X}$ is ergodic, with a finite and positive intensity given
by $\nu/(1-\mu)$.~Due to the branching structure, the number $S$
of offspring in a cluster follows the distribution
$(\ref{distrs})$.~Finally, we note that the classical Hawkes
process considered in the previous sections corresponds to the
special case where $d=1$ and $h(t)=0$ for $t\leq 0$.

A LDP for spatial Hawkes processes can be obtained by Theorem
\ref{scalarLDPspatPC}.~The precise statement is as Theorem
\ref{scalarLDPspatPC} with $(\ref{stability2})$ and

\[
\int_{\R^d}\|\xi\|h(\xi)\,\mathrm{d}\xi<\infty,
\]

\noindent in place of $(\ref{eq:hypPC})$ and $(\ref{eq:radCoSp})$,
moreover the rate function is $\Lambda^{*}(\cdot)$ defined by
$(\ref{rf})$.~Here the symbol $\|\cdot\|$ denotes the Euclidean
norm.

Similarly, the asymptotic behavior of the void probability
function and the empty space function of spatial Hawkes processes
can be obtained as immediate consequences of Proposition
\ref{cor:void} and Corollary \ref{cor:empty}, respectively.~The
precise statements are as Proposition \ref{cor:void} and Corollary
\ref{cor:empty}, with conditions $(\ref{stability2})$ and

\[
\int_{\R^d}\|\xi\|^d h(\xi)\,\mathrm{d}\xi<\infty
\]

\noindent in place of $\mathrm{E}[L^d]<\infty$.

\section{Extensions and open problems}
\label{sec:ext} In this paper we studied large deviations of
Poisson cluster processes.~Applications of these results to
insurance and queueing models are presently under investigation by
the authors.

The definition of Hawkes process extends immediately to the case
of random fertility rate $h(\cdot,Z)$, where $Z_{k}$'s are iid
unpredictable marks associated to the points $X_k$ (see Daley and
Vere-Jones (2003) for the definition of unpredictable marks, and
Br{\'e}maud, Massouli{\'e} and Ridolfi (2005) for the construction
of Hawkes processes with random fertility rate specified by an
unpredictable mark).~Due to the form of the distribution of $S$ in
this case (see formula $(6)$ in M{\o}ller and Rasmussen (2005)) it
is not clear if the LDPs for Hawkes processes proved in this paper
are still valid for Hawkes processes with random fertility rate.

The generalization of our results to non-linear Hawkes processes
(Kerstan (1964); Br{\'e}maud and Massouli{\'e} (1996);
Massouli{\'e} (1998); Br{\'e}maud, Nappo and Torrisi (2002);
Torrisi (2002)) would be interesting.~However, since a non-linear
Hawkes process is not even a Poisson cluster process, a different
approach is needed.

\vspace{0.3cm}

\section*{Acknowledgements}
We thank Kamil Szczegot for reporting two mistakes and for a careful reading of a first draft of the paper. 
\section*{References}

\noindent Borovkov, A.A. (1967), Boundary values problems for
random walks and large deviations for function spaces, {\em Theory
Probab. Appl.} {\bf 12}, 575--595.

\vspace{0.3cm}

\noindent Br{\'e}maud, P. (1981), {\em Point Processes and
Queues}, Springer, New York.

\vspace{0.3cm}

\noindent Br{\'e}maud, P. and Massouli{\'e}, L. (1996), Stability
of nonlinear Hawkes processes, {\em Ann. Prob.} {\bf 24},
1563--1588.

\vspace{0.3cm}

\noindent Br{\'e}maud, P., Massouli{\'e}, L. and Ridolfi, A.
(2005), Power spectra of random spike fields and related
processes, {\em Adv. Appl. Prob.} {\bf 37}, 1116--1146.

\vspace{0.3cm}

\noindent Br{\'e}maud, P., Nappo, G. and Torrisi, G.L. (2002),
Rate of convergence to equilibrium of marked Hawkes processes,
{\em J. Appl. Prob.} {\bf 39}, 123--136.

\vspace{0.3cm}

\noindent Brix, A. and Chadoeuf, J. (2002), Spatio-temporal
modeling of weeds by shot-noise G Cox processes, {\em Biometrical
J.} {\bf 44}, 83--99.

\vspace{0.3cm}

\noindent Chavez-Demoulin, V., Davison, A.C. and Mc Neil, A.J.
(2005), Estimating value-at-risk: a point process approach, {\em
Quantitative Finance} {\bf 5}, 227--234.

\vspace{0.3cm}

\noindent Daley, D.J. and Vere-Jones, D. (2003), {\em An
Introduction to the Theory of Point Processes (2nd edition)},
Springer, New York.

\vspace{0.3cm}

\noindent de Acosta, A. (1994), Large deviations for vector valued
L\'evy processes, {\em Stochastic Process. Appl.} {\bf 51},
75--115.

\vspace{0.3cm}

\noindent Dembo, A. and Zeitouni, O. (1998), {\em Large Deviations
Techniques and Applications (2nd edition)}, Springer, New York.

\vspace{0.3cm}

\noindent Ganesh, A., Macci, C. and Torrisi, G.L. (2005), Sample
path large deviations principles for Poisson shot noise processes,
and applications, {\em Electron. J. Probab.} {\bf 10}, 1026--1043.

\vspace{0.3cm}

\noindent Gusto, G. and Schbath, S. (2005), F.A.D.O.: a
statistical method to detect favored or avoided distances between
occurrences of motifs using the Hawkes model, {\em submitted}.

\vspace{0.3cm}

\noindent Hawkes, A.G. (1971a), Spectra of some self-exciting and
mutually exciting point processes, {\em Biometrika} {\bf 58},
83--90.

\vspace{0.3cm}

\noindent Hawkes, A.G. (1971b), Point spectra of some mutually
exciting point processes, {\em J. Roy. Statist. Soc. Ser. B} {\bf
33}, 438--443.

\vspace{0.3cm}

\noindent Hawkes, A.G. and Oakes, D. (1974), A cluster
representation of a self-exciting process, {\em J. Appl. Prob.}
{\bf 11}, 493--503.

\vspace{0.3cm}

\noindent Jagers, P. (1975), {\em Branching Processes with
Biological Applications}, John Wiley, London.

\vspace{0.3cm}

\noindent Jonnson, D.H. (1996), Point Process Models of
Single-Neuron Discharges, {\em J. Computational  Neuroscience}
{\bf 3}, 275--299.

\vspace{0.3cm}

\noindent Kerstan, J. (1964), Teilprozesse Poissonscher Prozesse.
{\em Transactions of the Third Prague Conference on Information
Theory, Statistical Decision Functions, Random Processes},
377--403.

\vspace{0.3cm}

\noindent Massouli{\'e}, L. (1998), Stability results for a
general class of interacting point processes dynamics, and
applications, {\em Stoch. Proc. Appl.} {\bf 75}, 1--30.

\vspace{0.3cm}

\noindent M{\o}ller, J. (2003), A comparison of spatial point
process models in epidemiological applications.~In {\em Highly
Structured Stochastic systems}, eds P.J. Green, N.L. Hjort and S.
Richardson, Oxford University Press, 264--268.

\vspace{0.3cm}

\noindent M{\o}ller, J. and Waagepetersen, R.S. (2004), {\em
Statistical Inference and Simulation for Spatial Point Processes},
Chapman and Hall, Boca Raton.

\vspace{0.3cm}

\noindent M{\o}ller, J. and Rasmussen, J.G. (2005), Perfect
simulation of Hawkes processes, {\em Adv. Appl. Prob.} {\bf 37},
629--646.

\vspace{0.3cm}

\noindent M{\o}ller, J. and Torrisi, G.L. (2007), The pair
correlation function of spatial Hawkes processes, {\em Statistics
and Probability Letters}, to appear

\vspace{0.3cm}

\noindent Neyman, J. and Scott, E.L. (1958), Statistical approach
to problems of cosmology, {\em J.R. Statist. Soc. B} {\bf 20},
1--43.

\vspace{0.3cm}

\noindent Ogata, Y. and Akaike, H. (1982), On Linear Intensity
Models for Mixed Doubly Stochastic Poisson and Self-Exciting Point
Processes, {\em J. R. Statist. Soc. B} {\bf 44}, 102--107.

\vspace{0.3cm}

\noindent Ogata, Y. (1988), Statistical models for earthquake
occurrences and residual analysis for point processes, {\em J.
Amer. Statist. Assoc.} {\bf 83}, 9--27.

\vspace{0.3cm}

\noindent Ogata, Y. (1998), Space-time point process model for
earthquake occurrences, {\em Ann. Inst. Statist. Math.} {\bf 50},
379--402.

\vspace{0.3cm}

\noindent Reynaud-Bouret, P. and Roy, E. (2007), Some non
asymptotic tail estimates for Hawkes processes, {\em Bull. Belg.
Math. Soc. Simon Stevin} {\bf 13}, 883--896.

\vspace{0.3cm}

\noindent Torrisi, G.L. (2002), A class of interacting marked
point processes: rate of convergence to equilibrium, {\em J. Appl.
Prob.} {\bf 39}, 137--161.

\vspace{0.3cm}

\noindent Vere-Jones, D. and Ozaki, T. (1982), Some Examples of
Statistical Estimation Applied to Earthquake Data, {\em Ann. Inst.
Statist. Math.} {\bf 34}, 189--207.

\end{document}